\documentclass[11pt]{article}
\usepackage{fullpage}

\usepackage{srcltx}
\usepackage{amssymb}
\usepackage{amsfonts}
\usepackage{graphicx}
\usepackage{bm}
\usepackage{theorem}
\usepackage{srcltx}
\usepackage{amsmath}
\usepackage{epsfig}
\usepackage{color}
\usepackage[colorlinks=true]{hyperref}

\def\Argmin{{\hbox{\rm Argmin}}}
\def\RIP{{\hbox{\rm RIP}}}

\def\cU{{\cal U}}

\def\cX{{\cal X}}
\def\cZ{{\cal Z}}
\def\cP{{\cal P}}

\def\cN{{\cal N}}
\def\sign{\mathop{\hbox{\rm sign}}}
\newtheorem{lemma}{Lemma}[section]
\newtheorem{proposition}{Proposition}[section]
\newtheorem{corollary}{Corollary}[section]

\newtheorem{algorithm}{Algorithm}%[section]

\def\bS{{\mathbf{S}}}

\def\bSG{{\mathbf{SG}}}
\def\bVSG{{\mathbf{VSG}}}
\def\bR{{\mathbf{R}}}
\def\Card{{\hbox{\rm Card}}}

\def\qed{$\Box$}

\def\Ker{{\hbox{\rm Ker}}}
\def\Opt{{\hbox{\rm Opt}}}
\def\Tr{{\hbox{\rm Tr}}}
\def\cU{{\cal U}}

\def\cC{{\cal C}}
\def\Conv{{\hbox{\rm Conv}}}

\def\e{\varepsilon}

\def\Argmax{{\hbox{\rm Argmax}}}
\newcommand{\be}{\begin{eqnarray}}
\newcommand{\ee}[1]{\label{#1}\end{eqnarray}}
\newcommand{\nn}{\nonumber \\}
\newcommand{\ese}{\end{eqnarray*}}
\newcommand{\bse}{\begin{eqnarray*}}
\newcommand{\rf}[1]{~(\ref{#1})}

\begin{document}

\title{Verifiable conditions of $\ell_1$-recovery for sparse signals with sign restrictions\thanks{Research of the second and the third authors was supported by the Office of Naval Research grant \# N000140811104.}}
%\subtitle{Do you have a subtitle?\\ If so, write it here}

%\titlerunning{Short form of title}        % if too long for running head

\author{Anatoli Juditsky         \and
        Fatma K{\i}l{\i}n\c{c} Karzan \and
        Arkadi Nemirovski
}

%\authorrunning{Short form of author list} % if too long for running head

\author{A. Juditsky \\
              LJK, Universit\'e J. Fourier, B.P. 53,
              38041 Grenoble Cedex 9, France \\
%              Tel.: +123-45-678910\\
%              Fax: +123-45-678910\\
              {\tt Anatoli.Juditsky@imag.fr}           %  \\
%             \emph{Present address:} of F. Author  %  if needed
           \and
           F. K{\i}l{\i}n\c{c} Karzan\\
              Georgia Institute of Technology,
              Atlanta, Georgia 30332, USA \\
              {\tt fkilinc@isye.gatech.edu}           %  \\
           \and
           A. Nemirovski\\
              Georgia Institute of Technology,
              Atlanta, Georgia 30332, USA \\
              {\tt nemirovs@isye.gatech.edu}           %  \\
}

%\date{Received: date / Accepted: date}

\maketitle

\begin{abstract}
We propose necessary and sufficient conditions for a sensing matrix to be ``$s$-semigood'' -- to allow for exact $\ell_1$-recovery of sparse signals with at most $s$ nonzero entries under sign restrictions on part of the entries.
We express error bounds for imperfect $\ell_1$-recovery
in terms of the characteristics underlying these conditions. These characteristics, although difficult to evaluate, lead to verifiable sufficient conditions for exact sparse $\ell_1$-recovery and thus efficiently computable upper bounds on those $s$ for which a given sensing matrix is $s$-semigood. We examine the properties of proposed verifiable sufficient conditions, describe their limits of performance and provide numerical examples comparing them with other verifiable conditions from the literature.
\end{abstract}

\section{Introduction}

Assessing a sparse signal from an observation has been one of the main research areas in Compressed Sensing and sparse signal recovery. In practice, {\em a priori} information about the signal to be recovered often exists and will be beneficial if taken into account in the recovery procedure. In this paper, we suppose that the {\em a priori} information about a {\em sparse} signal $w\in \bR^n$ amounts to the {\em sign restrictions}, and is given as the subsets $P_+$ and $P_-$ of $\{1,...,n\}$,  $P_+\cap P_-=\emptyset$, such that $w_i\ge 0$ for $i\in P_+$ and $w_i\le 0$ for $i\in P_-$. Therefore we address the following recovery problem: given an observation $y\in\bR^m$,
\be
y=Aw+e,\;\;\;
\ee{obs}
where $A\in \bR^{m\times n}$ (in this context $m< n$) is a given matrix, $e\in \bR^m$ is the observation error,
assess a {\em sparse} signal $w\in \bR^n$ satisfying {\em sign restrictions}.

A celebrated solution to the problem is given by the $\ell_1$-recovery, which amounts to taking, as an estimate of $w$, an optimal solution $\widehat{w}$ to the optimization problem
\be
\widehat{w}\in \Argmin_{x} \left\{\|x\|_1: \;\;\|Ax-y\|\le \e, \;x_i\ge 0\;\forall i \in P_+, ~x_i\le 0\;\forall i \in P_-\right\}
\ee{l1inex}
(here $\e$ is an {\em a priori} bound on the norm $\|e\|$ of the observation error, $\|\cdot\|$ being some norm on $\bR^m$).
When there are no sign restrictions (i.e. $P_+ = P_- = \emptyset$), we arrive at the estimator playing the central role in the Compressive Sensing theory. The central result here is that when signal $w$ is $s$-sparse (i.e., with at most $s$ nonzero entries) and the matrix $A$ possesses
a certain well-defined (although difficult to verify) property, then the $\ell_1$-recovery $\widehat{w}$ is close to $w$, provided the error bound $\e$ is small (for a comprehensive  survey see \cite{Candes_06} and references therein). Our goal here is to  propose efficiently verifiable sufficient conditions on $A$ which allow for similar `consistency'' results, with emphasis on the case where sign restrictions are present.

To outline our results and to position them with respect to what is already known, let us start with noiseless recovery (i.e., $\e=0$ and $y=Aw$). Here we are interested to answer the question:
\begin{quotation}
{\em Whether $A$ is such that whenever the true signal $w$ in {\rm \rf{obs}} is $s$-sparse and satisfies the sign constraints $w_i\geq0,$ $i\in P_+$, $w_i\leq0$, $i\in P_-$, the $\ell_1$-recovery
\be
\widehat{w}\in \Argmin_{x} \left\{\|x\|_1: \;\;Ax=y,\;x_i\ge 0\;\forall i \in P_+, ~x_i\le 0\;\forall i \in P_-\right\}
\ee{l1ex}
recovers $w$ exactly.
}\end{quotation}
If the answer is positive, we say that $A$  is {\sl $s$-semigood}\footnote{We use the term ``$s$-semigoodness'' to comply with the terminology of the  companion paper \cite{JNCS}, where  we used the name  {\sl $s$-goodness} to indicate that $\ell_1$-recovery as in \rf{l1ex} {\em without} the sign restrictions is exact.}.

The theory of Compressive Sensing provides several sufficient/necessary and sufficient conditions for the $\ell_1$-recovery to be exact. For example, when no sign constraints are imposed on $w$,
Donoho and Huo \cite{donhuo} prove that
$A$ is $s$-good if  for any set $I\subset \{1,...,n\}$ of cardinality $\leq s$ it holds
\be
\sum_{i\in I} |z_i|<\sum_{i\not\in I} |z_i|\;\;\mbox{for any}\;\;z\in \Ker A.
\ee{zhang1}
This condition has been extensively investigated. Its necessity has been established in \cite{Donoho_Elad_03}; it has been discussed in \cite{Zhang1,Zhang08} (under the name of {\em strict  $s$-balancedness}), where its link to the geometric necessary and sufficient condition of $s$-goodness from \cite{DonTan2} has been discussed. In \cite{Cohen_Dahmen_DeVore_06},  this condition has was  also related to the sufficient condition  (``{\em Null Space Property}'') for successful combinatorial recovery.

The first characterization of $s$-semigoodness for the case when $w$ is nonnegative (i.e. $P_+=\{1,...,n\}$) was proposed in
the founding paper of Donoho and Tanner \cite{DonTan1}
in terms of neighboring properties of the polytope $AS$, $S$ being the standard simplex $S=\{x\in \bR^n:\;x\ge 0, \;\sum_i x_i\le 1\}$. This paper contains also several important examples of $m\times n$ matrices which are $\lfloor {m\over 2}\rfloor$-semigood (here $\lfloor a\rfloor$ stands for the integer part of $a$) and demonstrates that various types of randomly generated matrices possess this property with overwhelming probability.
Extending the results from Donoho and Huo \cite{donhuo}, an equivalent characterization of $s$-semigoodness  has been provided in the nonnegative case by Zhang in \cite{Zhangnn,Zhang08}, where it is shown that $A$ is $s$-semigood if and only if the kernel of $A$, $\Ker A$, is {\sl strictly half $s$-balanced}, meaning that for any set $I\subset \{1,...,n\}$ of cardinality $\leq s$ it holds
\be
\sum_{i\in I} z_i<\sum_{i\not \in {I}} |z_i|\;\;\mbox{for any}\;\;z\in \Ker A\;\mbox{such that}\; z_i\le 0,\;\mbox{for all}\;i\not\in {I}.
\ee{zhangnn}

It should be mentioned that the necessary and sufficient conditions for $s$-semigoodness from  \rf{zhang1}, \rf{zhangnn} and \cite{DonTan1,DonTan2} share a common drawback -- they seemingly cannot be verified in a computationally efficient way. To the best of our knowledge, the only {\sl efficiently verifiable} conditions for $s$-semigoodness offered by the existing Compressive Sensing theory are the {\sl sufficient} conditions based on the {\em mutual incoherence}
\begin{equation}\label{muin}
\mu(A)=\max_{i\neq j}{|A_i^TA_j|\over A_i^TA_i}
\end{equation}
where $A_i$ are columns of $A$ (assumed to be nonzero).
Clearly, the mutual incoherence can be easily computed even for large matrices. Unfortunately, it turns out that that the estimates of ``level of (semi)goodness'' of a sensing matrix based on mutual incoherence usually are too conservative, in particular, they are provably dominated by the verifiable Linear Programming (LP) based sufficient conditions for $s$-goodness proposed in the companion paper \cite{JNCS} and based on characterization of $s$-goodness given in \rf{zhang1}. Another verifiable sufficient condition for $s$-goodness, which uses the Semidefinite Programming (SDP) relaxation,  has been recently proposed in \cite{AspGaou}.

The contributions of this paper, which follow the approach developed in \cite{JNCS},  are as follows.
\begin{enumerate}
\item Taking existing characterizations of (semi)goodness \rf{zhang1}, \rf{zhangnn} as a starting point, we develop in Section \ref{sec:2}, several equivalent necessary and sufficient conditions for $s$-semigoodness of a matrix $A$ in the case of general-type sign restrictions. Then in Section \ref{sec:3}, we establish error bounds for inexact $\ell_1$-recovery (noisy observation \rf{obs}, imprecise optimization in \rf{l1inex}, nearly-sparse true signals); these bounds are expressed  in the same terms as the necessary and sufficient conditions for $s$-semigoodness from Section \ref{sec:2}. These bounds can be seen as an extension to the sign restricted case of bounds of Section 3 in \cite{JNCS} and as a special case of the bounds provided in Theorem 4.1 of \cite{Zhang08}. To the best of our knowledge, these bounds that incorporate sign information of the signal are new.
\item The major goal of this paper is to use the LP relaxation techniques from \cite{JNCS} to derive novel {\em efficiently verifiable} sufficient conditions for $s$-semigoodness. These conditions allow  one to build, in a computationally efficient fashion, lower bounds on the ``level of $s$-semigoodness'' of a given matrix $A$, that is, on the largest $s=s_*(A)$ for which $A$ is $s$-semigood with respect to given $P_\pm$. Some properties of these verifiable conditions, same as limits of their performance, are studied in Sections \ref{sec:4}, \ref{sec:limits}, where we provide also a computationally efficient scheme for upper bounding of $s_*(A)$. In Section \ref{sec:sdp}, we develop another efficiently computable lower bound for $s_*(A)$ by applying the SDP relaxation, similar to the approach developed in \cite{AspGaou} for the ``unsigned'' case $P_\pm=\emptyset$. In Section \ref{sec:numerical} we report on numerical experiments aimed at comparing the ``power'' of our LP-based sufficient conditions for $s$-semigoodness,
their ``unsigned'' prototypes from \cite{JNCS}, and conditions based on mutual incoherence. We show that incorporating the sign information can improve the bounds on the level of $s$-semigoodness, and that the bounds based on LP relaxations clearly outperform the bounds based on mutual incoherence.
\item It turns out that our verifiable sufficient conditions for $s$-semigoodness can be expressed in terms of specific properties of the {\sl linear} recovery $\widehat{w}^{\hbox{\scriptsize lin}}=Y^Ty$ associated with an appropriate $m\times n$ matrix $Y$. In Section \ref{sec:mp}, we propose and justify a new non-Euclidean {\em Matching Pursuit} algorithm associated with this linear recovery.
\end{enumerate}
\section{Necessary and sufficient conditions for $s$-semigoodness}
\label{sec:2}
Let $A$ be an $m\times n$ matrix, let $s$, $1\leq s\leq m$, be an integer, and let $P_+, ~P_-$ and $P_n$ be a partition of $\{1, \ldots, n\}$ into three non-overlapping subsets. We say that $A$ is {\em $s$-semigood}, if for every vector $w$ with at most $s$ nonzero entries satisfying $w_i\geq0$ for $i\in P_+$, and $w_i\leq0$ for $i\in P_-$, $w$ is the unique optimal solution to the problem
\be
\Opt = \min_z \left\{ \|z\|_1: Az=Aw, ~z_i \ge 0 ~\forall i \in P_+, ~z_i \le 0~ \forall i \in P_- \right\}.
\ee{OptPrb1}
\par
Our primary goals are to find necessary and sufficient and {\sl verifiable} sufficient conditions for $A$ to be $s$-semigood.
\par
Note that without loss of generality we may assume $P_-=\emptyset$. Indeed, by replacing the partition $P_+,\,P_-,P_n$ with the partition
$\overline{P}_+=P_+\cup P_-,\overline{P}_-=\emptyset,\overline{P}_n=P_n$ and matrix $A$ -- with the matrix $\overline{A}$ obtained from $A$ by multiplying the columns with indices $i\in P_-$ by $-1$, $s$-semigoodness of $A$ with respect to the original sign restrictions given by $P_\pm,P_n$ is equivalent to the $s$-semigoodness of the new matrix $\overline{A}$ with respect to the new sign restrictions.
By this reason, {\sl we assume from now on that $P_-=\emptyset.$} Besides this, we assume without loss of generality that $P_+=\{1,...,p\}$ and $P_n=\{p+1,...,n\}$ for some $p$. From now on, we denote by $\cP_n$ the set of all signals satisfying the sign restrictions:
$$
\cP_n=\{w\in\bR^n: w_i\geq0\,~\forall i\in P_+\}.
$$
\par
Note that since $P_-=\emptyset$, \rf{OptPrb1} simplifies to
\be
\Opt = \min_z \left\{ \|z\|_1: Az=Aw, ~z_i \ge 0 ~\forall i \in P_+\right\}.
\ee{OptPrb}
\par
Let us fix a norm $\|\cdot\|$ on $\bR^n$, and let $\|\cdot\|_*$ be the conjugate norm.
\begin{proposition}\label{prop1} Let $m,n,s$ and $P_+$ be given. The following six conditions on an $m\times n$ matrix $A$ are equivalent to each other:
\par
{\rm (i)} $A$ is $s$-semigood;
\par
{\rm (ii)} For every subset $J$ of $\{1,...,n\}$ with $\Card(J)\leq s$, and  any $x\in \Ker A\backslash \{0\}$ such that $x_i\leq0$ for all $i\in P_+\setminus J$ one has
\[
\sum_{i\in J\cap P_+}x_i + \sum_{i\in J\cap P_n}|x_i| < \sum_{i\not\in J}|x_i|.
\]
\par
{\rm (iii)} There exists $\xi\in(0,1)$ such that
for every subset $J$ of $\{1,...,n\}$ with $\Card(J)\leq s$ and any $x\in\Ker A$ such that
$x_i\leq0$ for all $i \in P_+ \setminus J$
one has
\[\sum_{i\in J\cap P_+}x_i + \sum_{i\in J\cap P_n}|x_i| \leq \xi\sum_{i\not\in J}|x_i|.
\]
\par
{\rm (iv)} There exist $\xi\in(0,1)$ and  $\theta\in[1,\infty)$ such that $A$ satisfies the condition $\bSG_{s}(\xi,\theta)$ as follows: \\
for every $x\in\Ker A$ and every subset $J$ of $\{1,...,n\}$ with $\Card(J)\leq s$, one has
\[
\sum_{i\in J\cap P_+}x_i + \sum_{i\in J\cap P_n}|x_i| \leq \xi \left( \sum_{i\in P_n\setminus J}|x_i|  + \sum_{i\in P_+ \setminus J}\psi(x_i) \right), \;\;\psi(t)=\max[-t,\theta t],
\]
or, equivalently: for all $x\in \Ker A$, $\Theta(x) \leq \xi\Psi(x)$ where
\begin{equation}\label{eq444}
\begin{array}{l}
\Theta(x):=\max\limits_{{J\subset \{1,...,n\},\atop\hbox{\scriptsize\rm Card}(J)\leq s}}\left[\sum_{i\in J\cap P_+}\max[(1-\xi)x_i,(1+\theta\xi)x_i]+\sum_{i\in J\cap P_n}(1+\xi)|x_i|\right]\\
\Psi(x):=\sum_{i\in P_+}\max[-x_i,\theta x_i]+\sum_{i\in P_n} |x_i|\\
\end{array}
\end{equation}
\par
{\rm (v)} There exist $\xi\in(0,1)$, $\theta\in[1,\infty)$ and $\beta\in[0,\infty)$ such that $A$ satisfies the condition $\bSG_{s,\beta}(\xi,\theta)$ as follows: \\
for every $x\in\bR^n$ and every subset $J$ of $\{1,...,n\}$ with $\Card(J)\leq s$, one has
\[
\sum_{i\in J\cap P_+}x_i + \sum_{i\in J\cap P_n}|x_i| \leq \beta\|Ax\|+\xi\left( \sum_{i\in P_n\setminus J}|x_i|  + \sum_{i\in P_+ \setminus J}\psi(x_i) \right), \;\;\psi(t)=\max[-t,\theta t].
\]
\par
{\rm (vi)} There exist $\xi\in(0,1)$ and  $\beta\in[0,\infty)$ such that $A$ satisfies the condition $\bSG_{s,\beta}(\xi)$ as follows: \\
for every $J\subset\{1,...,n\}$ with $\Card(J)\leq s$ and any $x\in\bR^n$ such that $x_i\leq0$ for all $i \in P_+ \setminus J$, one has
\[
\sum_{i\in J\cap P_+}x_i + \sum_{i\in J\cap P_n}|x_i| \leq \beta\|Ax\|+\xi\sum_{i\not\in J}|x_i|.
\]
\end{proposition}

We provide the proof of Proposition \ref{prop1} in Appendix \ref{App:Prop1}.

As we have already mentioned in Introduction, when $P_n = \emptyset$ or $P_+=\emptyset$, the characterizations (i)--{(iv)} of $s$-semigoodness are not completely new. For instance, when $P_n=\emptyset$, a necessary and sufficient condition for $s$-semigoodness of $A$ in the form {\rm (ii)} has been established in \cite{Zhangnn} (compare {\rm (ii)} to the definition \rf{zhangnn} of half $s$-balancedness of $\Ker A$).
On the other hand, the equivalent formulation of this characterization in terms of conditions $\bSG_{s,\beta}(\xi,\theta)$ and $\bSG_{s,\beta}(\xi)$ seems to be new. We are about to demonstrate that the latter two conditions allow to control the error of $\ell_1$-recovery in the case when the vector $w\in \bR^n$ is not $s$-sparse and the problem
 \eqref{OptPrb} is not solved to exact optimality.

\section{Error bounds for imperfect $\ell_1$-recovery}
\label{sec:3}
We have seen that the conditions provided in Proposition \ref{prop1} are responsible for $s$-semigoodness of a sensing matrix $A$, that is, for the exactness of $\ell_1$-recovery in the ``ideal case'' when the true signal $w$ is $s$-sparse, there is no observation error, and the optimization problem \rf{OptPrb} is solved to exact optimality. Below we demonstrate that these conditions control also the error of $\ell_1$-recovery in the case when the signal $w\in\cP_n$ is not exactly $s$-sparse, there is observation noise and problem \rf{OptPrb} is not solved to exact optimality. The corresponding error bound (cf \cite[Proposition 3.1, Theorem 3.1]{JNCS}) is as follows:

\begin{proposition}\label{IRB}  Let  $w\in\cP_n$ be such that $\|w-w^s\|_1\leq\mu$, where $w^s$ is the vector obtained from $w$ by replacing all but the $s$ largest in magnitude entries in $w$ with zeros, let $y$ be such that $\|Aw-y\|\leq\e$, and let, finally, $x$ be an approximate solution to the optimization problem
\begin{equation}\label{problem!}
\Opt=\min\limits_z\left\{\|z\|_1: \|Az-y\|\leq\e, ~z_i\geq0 ~\forall i \in P_+\right\}.%\eqno{(!)}
\end{equation}
such that $\|x\|_1\leq\Opt+\nu$ and $\|Ax-y\|\leq\delta$.
\begin{enumerate}
\item If $A$ satisfies the condition $\bSG_{s,\beta}(\xi,\theta)$ with some $\xi\in(0,1)$, $\beta\in[0,\infty)$ and $\theta\in [1,\infty)$, then
\be
\|x-w\|_1\leq
\frac{1+\xi}{1-\xi} \nu + \frac{2(1+\xi\theta)}{1-\xi}\mu + \frac{2\beta}{1-\xi} (\e+\delta).
\ee{eq256}
\item If $A$ satisfies the condition  $\bSG_{s,\beta}(\xi)$ with some $\xi\in(0,1)$,  $\beta\in[0,\infty)$, then
    \be
\|x-w\|_1\leq\frac{1+\xi}{1-\xi} \nu + \frac{{2(1+\beta\alpha)}}{1-\xi}\mu + \frac{2\beta}{1-\xi} (\e+\delta).
\ee{eq256'}
where $\alpha$ stands for the maximum of $\|\cdot\|$-norms of the columns in $A$.
\end{enumerate}
\end{proposition}
For proof, see Appendix \ref{App:IRB}.

\section{Verifiable conditions for $s$-semigoodness}
\label{sec:4}
We are about to demonstrate that condition  ${ \bSG_{s,\beta}(\xi,\theta)}$ from Proposition \ref{prop1} leads to efficiently computable lower and upper bounds on the level of $s$-semigoodness.

\subsection{Verifiable sufficient conditions for $s$-semigoodness by Linear Programming}
\label{LowerBounds}
Let
\[
\cU_s=\left\{u\in\bR^n:~\|u\|_1\leq s,~ \|u\|_\infty \leq1\right\},
\]
so that $\cU_s$ is the convex hull of all $\{-1,0,1\}$ vectors with at most $s$ nonzero entries, and for $x\in \bR^n$, let $\|x\|_{s,1}$ be the sum of the $s$ largest magnitudes of entries in $x$, or, equivalently,
$$
\|x\|_{s,1}=\max_{u\in \cU_s}u^Tx.
$$
Let
$$
(D_\theta[x])_i=\left\{\begin{array}{ll}[1+\theta\xi]\max[x_i,0],&i\in P_+\\
(1+\xi)|x_i|,&i\not\in P_+\\
\end{array}\right.,\qquad \Phi(x)=\|D_\theta[x]\|_{s,1}.
$$
Suppose $\xi\in[0,1)$, $\theta\in[1,\infty)$ and $\rho,\sigma\in[0,\infty)$ are given. Consider the following condition on an $m\times n$ matrix $A$:
\begin{quote}
$\bVSG_s(\xi,\theta,\rho,\sigma)$: {\sl There exist $m\times n$ matrix $Y=[y_1,...,y_n]$ and a vector $v\in\bR^m$ such that
\begin{equation}
\begin{array}{rclr}
\Phi_s(-C_i[Y,A])+(A^Tv)_i &\leq& \xi,\,\,\,1\leq i\leq n&(a)\\
\Phi_s(C_i[Y,A])-(A^Tv)_i &\leq& \xi,\,\,\,i\not\in P_+&(b)\\
\Phi_s(C_i[Y,A])-(A^Tv)_i &\leq& \theta \xi,\,\,i\in P_+&(c)\\
\|y_i\|_* &\leq& \sigma,\,\,\,1\leq i\leq n&(d)\\
\|v\|_* &\leq& \rho&(e)
\end{array}\label{condition}
\end{equation}
where $C_i[Y,A]$ is the $i$-th column of the matrix $I-Y^TA$.}
\end{quote}
\par
Observe that this condition is verifiable, since (\ref{condition}) is a system of explicit convex constraints on $Y$ and $v$.
\begin{proposition}\label{propmain} Let $A$ satisfy $\bVSG_s(\xi,\theta,\rho,\sigma)$ with some $\xi\in[0,1)$, $\theta\in[1,\infty)$, and $\rho,\sigma\in[0,\infty)$. Then $A$ satisfies $\bSG_{s,\beta}(\xi,\theta)$ with
\begin{equation}\label{beta}
\beta=\rho+\sigma\max\limits_{k_+,k_n}\left\{k_+(1+\theta\xi)+k_n(1+\xi):\begin{array}{l}
0\leq k_+\leq \Card(P_+)\\
0\leq k_n\leq\Card(P_n)\\
k_++k_n\leq s\\
\end{array}\right\}\leq \rho+\sigma s(1+\theta\xi).
\end{equation}
In particular, $A$ is $s$-semigood.
\end{proposition}
For proof, see Appendix \ref{App:Propmain}.

Some comments are in order.
\paragraph{Origin of the condition $\bSG_{s,\beta}(\xi,\theta)$.}
The condition $\bVSG_s(\xi,\theta,\rho,\sigma)$ is yielded by a simple and general construction, and we believe it makes sense to present this construction in its general form. The essence of the matter is in building a verifiable sufficient condition for the validity of (\ref{eq444}), see Proposition \ref{prop1}.iv. By positive homogeneity of degree 1  of the convex functions $\Theta,\Psi$ participating in (\ref{eq444}), the latter condition is equivalent to
\begin{equation}\label{equivto}
\Opt:=\max\limits_x\left\{\Theta(x):Ax=0,x\in X\right\}\leq\xi, \,\,\,\,X=\{x:\Psi(x)\leq1\}.
\end{equation}
A verifiable sufficient condition for (\ref{equivto}) is basically the same as an efficiently computable upper bound for $\Opt$; the sufficient condition for the validity of (\ref{equivto}) associated with such a bound merely states that the bound is $\leq\xi$. Now observe that from the origin of $\Psi$ (see (\ref{eq444})) it is clear that $X$ has a moderate number, $N$, of readily available extreme points $x^1,...,x^N$ (in the case of (\ref{eq444}), $N=2n$), so that the only difficulty in computing $\Opt$ exactly comes from linear constraints $Ax=0$. The standard way to circumvent this difficulty and to efficiently bound $\Opt$ from above  is to use the Lagrange relaxation: for any $v\in \bR^m$,
\bse
\Opt&=&\max\limits_{x\in X} \left\{\Theta(x)+v^TAx:Ax=0,x\in X\right\}\\
&\leq &\max\limits_x\left\{\Theta(x)+v^TAx:x\in X\right\}=
\max\limits_{1\leq i\leq N} [\Theta(x^i)+v^TAx^i],
\ese
and hence the efficiently computable {\sl Lagrange relaxation bound} $\inf_v\max_{1\leq i\leq N} [\Theta(x^i)+v^TAx^i]$ is an upper bound on $\Opt$. Unfortunately, in our situation this bound can be very poor; e.g., when $X$ is symmetric with respect to the origin and $\Theta$ is even (as it happens in (\ref{eq444}) when $P_+=\emptyset$), it is immediately seen that the bound becomes the trivial bound $\Opt\leq\max_{x\in X} \Theta(x)=\max_i\Theta(x^i)$. In order to strengthen the relaxation, we pass to the Fenchel-type representation of $\Theta$
$$
\Theta(x)=\sup\limits_u\left[[Pu+q]^Tx-\Theta_*(u)\right]
$$
with a proper convex function $\Theta_*$; such a representation, even with $Pu+p\equiv u$, exists whenever $\Theta$ is a proper convex function (and can be easily found for $\Theta$ we are interested in). We now have for any $Y\in \bR^{m\times n}$, $v\in\bR^m$,
\bse
\Opt&=&\max\limits_x\left\{\Theta(x):Ax=0,x\in X\right\}=\sup\limits_{x,u}\left\{[Pu+p]^Tx-\Theta_*(u):Ax=0,x\in X\right\}\\
&=&
\sup\limits_{x,u}\left\{[Pu+p]^T[x-Y^TAx]+v^TAx-\Theta_*(u): Ax=0,x\in X\right\}\\
&\leq& \sup\limits_{x,u}\left\{[Pu+p]^T[x-Y^TAx]+v^TAx-\Theta_*(u): x\in X\right\}\\
& =&\max\limits_{1\leq i\leq N}\underbrace{\sup\limits_u\left\{[Pu+p]^T[x^i-Y^TAx^i]+v^TAx^i-\Theta_*(u)\right\}}_{:=\Theta_i(Y,v)},
\ese
so that the condition
\begin{equation}\label{thecondition}
\exists(Y\in\bR^{m\times n},v\in\bR^m): \Theta_i(Y,{ v})\leq\xi,\,\,1\leq i\leq N,
\end{equation}
is sufficient for the validity of (\ref{equivto}). Note that the functions $\Theta_i$, by their origin, are convex, so that the condition (\ref{thecondition}) is efficiently verifiable, provided that $\Theta_i(\cdot)$ are efficiently computable.
\par
In the case we are interested in, the extreme points of $X$ are the $2n$ vectors $-e_i$ for $1\leq i\leq n$, $e_i$ for $i\in P_n$, and $\theta^{-1}e_i$ for $i\in P_+$, where $e_i$ is the $i$-th basic orth. Implementing the outlined bounding scheme and adding additional restrictions (\ref{condition}.$d$,$e$) to get a control over $\beta$, we arrive at (\ref{condition}). It should be stressed that the outlined scheme can be applied to bounding from above the optimal value of a whatever problem of the form (\ref{equivto}) with a convex polytope $X$ and a proper convex objective $\Theta$; all what matters is that $X$ is given as $\Conv\{x^1,...,x^N\}$ and $\Theta$ is efficiently computable. Note also that when $X$ is a polytope given by list of $M$ linear inequalities, we can efficiently represent it as the intersection of $M$-dimensional standard simplex and an affine plane, so that the outlined scheme is applicable to a whatever problem of maximizing an efficiently computable proper convex function under a (finite) system of linear inequality and equality constraints.
\paragraph{Effect of increasing $\beta,\theta,\xi$.} The condition $\bSG_{s,\beta}(\xi,\theta)$ appearing in Proposition \ref{prop1}.v clearly is ``monotone'' in the parameters $\beta,\theta,\xi$: whenever $A$
satisfies this condition and $\beta'\geq\beta$, $\theta'\geq\theta$ and $\xi'\geq\xi$, $A$ satisfies the condition $\bSG_{s,\beta'}(\xi',\theta')$ as well. Proposition \ref{propmain}
offers a verifiable sufficient condition for the validity of $\bSG_{s,\beta}(\xi,\theta)$, specifically,
\begin{quote} $\bVSG^*_{s,\beta}(\xi,\theta)$: $\exists Y,v$ $\rho$, $\sigma$ satisfying (\ref{condition}) and the relation ${ \rho+ \sigma s(1+\theta\xi)}\leq\beta$.
\end{quote}
A natural question is, whether
this verifiable condition possesses the same monotonicity properties as the ``target'' condition $\bSG_{s,\beta}(\xi,\theta)$. In the case of the affirmative answer, in order to
conclude that $A$ is $s$-semigood, we could check the validity of $\bVSG^*_{s,\beta}(\xi,\theta)$ for appropriately large values of $\beta,\theta$ and a close to one value of $\xi<1$; if the condition is satisfied, $A$ is $s$-semigood, and error bounds
from Proposition \ref{IRB} take place. Were the condition $\bVSG^*_{s,\beta}(\xi,\theta)$ ``not monotone,'' to justify the $s$-semigoodness of $A$ via this condition would require a problematic and time-consuming search in the space of parameters $\beta,\theta,\xi$. Fortunately, the condition $\bVSG^*_{s,\beta}(\xi,\theta)$ indeed is monotone:
\begin{proposition}\label{propnew} Let $A$ satisfy $\bVSG^*_{s,\beta}(\xi,\theta)$, and let $Y,v,\sigma,\rho$ be the corresponding certificate, that is, ${ \rho+ \sigma s(1+\theta\xi)}\leq\beta$ and $Y,v,\sigma,\rho$ satisfy {\rm (\ref{condition})}. Then $A$ satisfies $\bVSG^*_{s,\beta'}(\xi',\theta')$  whenever $\beta'\geq\beta$, $\theta'\geq\theta$ and $\xi'\in(\xi,1)$, the certificate being $(Y',v,\sigma,\rho)$, where the columns $Y_i^\prime$ of $Y'$ are multiplies of the columns $Y_i$ of $Y$, namely,
$$
Y_i^\prime=a_iY_i;\,\,\,[0,1]\ni a_i=\left\{\begin{array}{ll}(1+\xi\theta)/(1+\xi'\theta'),&i\in P_+\\
(1+\xi)/(1+\xi'),&i\in P_n\\
\end{array}\right.
$$
\end{proposition}
For proof, see {\color{blue} Online Supplement \ref{App:Monotonicity}}.

\paragraph{Relation to the sufficient condition for $s$-goodness from \cite{JNCS} and the Restricted Isometry Property.}
The verifiable sufficient condition for $s$-goodness from \cite{JNCS} requires from an $m\times n$ matrix $A$ the existence of $\gamma<1/2$ and $Y=[y_1,...,y_n]\in\bR^{m\times n}$ such that
$$
\|C_i[Y,A]\|_{s,1}\leq\gamma,\;\;\mbox{for all $1\leq i\leq n$},
$$
Setting $\theta=1$ and $\xi={\gamma\over 1-\gamma}$ (so that $\xi<1$ and $\gamma=\frac{\xi}{1+\xi}$) and taking into account that in the case of $\theta=1$ we have $\Phi_s(z)\leq (1+\xi)\|z\|_{s,1}$, the latter condition implies that
$$
\Phi_s(\pm C_i[Y,A]) \leq (1+\xi)\gamma = \xi, ~~\forall i,
$$
that is, it implies the validity of $\bVSG_s(\xi,1,0,\sigma)$, provided that $\sigma$ is large enough, specifically, $\sigma\geq\|y_i\|_*$ for all $i$.

As it was shown in the companion paper \cite{JNCS}, when $A$ satisfies the {\em Restricted Isometry Property} $\RIP(\delta,k)$ with parameters $\delta\in(0,1)$, $k>1$, the above sufficient condition for $s$-goodness is satisfied with $\gamma=1/3$ for $s$ as large as $O(1)(1-\delta)\sqrt{k}$; as a result, a $\RIP(\delta,k)$-matrix satisfies $\bVSG_s({1\over2},1,0,\sigma)$ provided that $\sigma$ is large enough and $s\leq O(1)(1-\delta)\sqrt{k}$. Since for large $m,n$, $m<n$, typical random matrices possess, with overwhelming probability, property $\RIP({1\over2},k)$ with $k$ as large  as $O(1)m/\ln(n/m)$, we see that our verifiable sufficient condition for $s$-semigoodness can certify the latter property for  $s$ as large as $O(1)\sqrt{m/\ln(n/m)}$, provided that the matrix in question is ``good enough''.

\subsection{Upper bounding the level of $s$-semigoodness} \label{sec:ub}\label{UpperBounds}
Here we address the issue of bounding from above the maximal $s=s_*(A)$ for which $A$ is $s$-semigood.
The construction to follow is motivated by item (iv) of Proposition \ref{prop1}. A necessary and sufficient condition for the $s$-semigoodness of $A$ is  the existence of $\xi<1$ and $\theta\ge 1$ such that for all $x\in\Ker A$ and any set $I$ of indices with $\Card(I)\le s$
\bse
\sum_{i\in I\cap P_+}\max[(1-\xi)x_i,(1+\theta\xi)x_i]+\sum_{i\in I\cap P_n}(1+\xi)|x_i|\leq \xi\Psi(x)
\ese
where
\be
\Psi(x)=\sum_{i\in P_+}\max[-x_i,\theta x_i]+\sum_{i\in P_n} |x_i|,
\ee{psi2}
or, equivalently,
\begin{quote}
(!) for every $x\in\Ker A$ and every vector $v$ with at most $s$ nonzero entries and nonzero entries $v_i$ belonging to $[1-\xi,1+\xi\theta]$ if $i \in P_+$ and belonging to $[-1-\xi,1+\xi]$ if $i \in P_n$,
one has
$$
v^Tx\leq \xi \Psi(x).
$$
\end{quote}
Observe that the convex hull of the vectors $v$ in question is exactly the set
$$
\cU^{\xi,\theta}=\left\{v\in\bR^n:
    \begin{array}{l}
    0\leq v_i\leq 1+\theta\xi, \,i\in P_+, ~|v_i|\leq 1+\xi,\,i\in P_n, \\
    ~\sum_{i\in P_+}{v_i\over1+\theta\xi}+\sum_{i\in P_n}{|v_i|\over1+\xi}\leq s
    \end{array}
\right\}.
$$
Recalling that $P_+=\{1,...,p\}$, setting $q=n-p=\Card(P_n)$ and
\be
\cU=\left\{u\in\bR^n:~\|u\|_1\le s, ~\|u\|_{\infty}\le 1, ~u_i \ge 0 \mbox{ for } i \in P_+ \right\}
\ee{uu}
we see that
\be
\cU^{\xi,\theta}=V^{\xi,\theta}\cU,~~ \textup{ where }~~ V^{\xi,\theta}=\left[\begin{array}{c|c}(1+\xi\theta)I_p&0\cr\hline 0&(1+\xi)I_q\cr\end{array}\right].
\ee{Vxite}
The condition (!) now reads
$$
\max_{v\in \cU^{\xi,\theta}} v^Tx\leq \xi \Psi(x)\;\mbox{ for all } x\in\Ker A.
$$
Setting $\cX=\{x\in \Ker A: ~\Psi(x) \le 1\}$ the latter condition, by homogeneity reason, is the same as
\be
\Opt=\Opt(\xi,\theta):=\max_{v,x}\left\{v^Tx:\; v\in\cU^{\xi,\theta},\; x\in\cX \right\}
\leq\xi;
\ee{ub1}
recall that $A$ is $s$-semigood if and only if there exist $\theta\geq1$ and $\xi<1$ such that \rf{ub1} takes place.

We can use \rf{ub1} in order to bound  $s_*(A)$ from above, as follows.
In order to certify that $s_*(A)<s$ for a given $s$ ($s$ is the input to our algorithm), we fix a large $\theta$ and a close to one $\xi<1$ (these are the parameters of the algorithm) and run the iterations
$$
u_0\in \cU^{\xi,\theta}\mapsto x_1\in\Argmax_{x\in \cX} u_0^Tx\mapsto u_1\in\Argmax_{u\in \cU^{\xi,\theta}}u^Tx_1\mapsto
%x_2\in\Argmax_{x\in \cX}u_1^Tx\mapsto
...
$$
initiating them by a picked at random vertex $u_0$ of $\cU^{\xi,\theta}$. Note that the quantities $u_i^Tx_i$, $i=1,2,...$  clearly form a nondecreasing sequence of lower bounds on $\Opt$. We terminate the outlined iterations when the progress in the bounds -- the difference $u_i^Tx_i-u_{i-1}^Tx_{i-1}$ -- falls below a given small threshold, and we run this process a predetermined number of times from different randomly chosen starting points. As a result, we get a set of lower bounds on $\Opt$ of the form $u^Tx$, where $u$ is a vertex of $\cU^{\xi,\theta}$ and $x\in \cX$. If our goal were merely to certify that \rf{sdp1} is not valid for given $s,\theta,\xi$, we could terminate this process at the first step, if any, when the current lower bound $u^Tx$ becomes $>\xi$ (cf. \cite[Section 4.1]{JNCS}). We, however, want to certify that $s>s_*(A)$, or, which is the same by Proposition \ref{prop1}.iv, that \rf{sdp1} fails to be true for {\sl all} $\theta$ and all $\xi<1$, and not only for those $\theta,\xi$ we have selected for our test. To overcome this difficulty, we accompany every step
$u\mapsto x\in \Argmax_{x\in\cX}u^Tx$ by an additional computation as follows. In our process, $u$ is an extreme point of $\cU^{\xi,\theta}$, that is, a point with $s_u\leq s$ nonzero entries, let the set of indices of these entries be $I$. Setting $\epsilon_i=\sign(u_i)$, we solve the following LP problem
$$
\max\limits_x\left\{\sum_{i\in I\cap P_+}x_i+\sum_{i\in I\cap P_n}\epsilon_ix_i:\left\{\begin{array}{l}x_i\leq0,~i\in P_+\backslash I\\
Ax=0\\
\sum_{i\not\in I}|x_i|\leq 1\\
\end{array}\right.\right\}.
$$
If the optimal value in this problem is $\geq1$, we terminate our test and claim that $A$ is not $s$-good; by Proposition \ref{prop1}.ii, this indeed is the case.
\par
As applied to a given input $s$, the outlined test either terminates with a valid claim ``$s>s_*(A)$'', or terminates with no conclusion at all, in which case we could pass to testing a larger value of $s$.

\section{Limits of performance of LP-based sufficient conditions for $s$-semi\-go\-od\-ness}
\label{sec:limits}
Unfortunately,  the condition in question, same as its predecessor from \cite{JNCS}, can{\em not} certify $s$-semi\-go\-od\-ness of an $m\times n$ matrix in the case of $s>O(1)\sqrt{m}$, unless the matrix is ``nearly square''. The precise statement is as follows (cf. \cite[Proposition 4.2]{JNCS}):
\begin{proposition}\label{limits} Let
\begin{equation}\label{eq119}
n>2(2\sqrt{2m}+1)^2
\end{equation}
and let $\xi<1,\theta\geq1,\sigma\geq0,\rho\geq0$, an integer $s$ and an $m\times n$ matrix $A$ be such that $A$ satisfies $\bVSG_s(\xi,\theta,\rho,\sigma)$. Then
\begin{equation}\label{eq191}
s\leq 2\sqrt{2m}+1.
\end{equation}
\end{proposition}

For proof, see Appendix \ref{App:Proplimits}.

The results from Proposition \ref{limits} show that our verifiable sufficient conditions can only certify $s$-semigoodness of an $m \times n$ matrix at a suboptimal rate of $s\le O(1)\sqrt{m}$, unless the matrix is ``nearly square''. In fact this verifiable bound can still give a very poor impression on the true largest $s=s_*(A)$ for which $A$ is $s$-semigood.
An instructive example in this direction is as follows. Consider the case of $P_+=\{1,...,n\}$, let $m=2d+1$ be odd, and let the rows of $A$ be comprised of the values of basic trigonometric polynomials
\[
p_0(\phi)\equiv 1, \;\; p_{2i-1}(\phi)=\cos(i\phi), \;\;p_{2i}(\phi)=\sin(i\phi), \;\;1\leq i\leq d,
\]
taken along the regular grid $\phi_j=2\pi j/n$, $0\leq j<n$, so that $A_{ij}=p_i(\phi_j)$, $0\leq i<m$, $0\leq j<n$ (we enumerate rows and columns starting with 0 rather than with 1). It is well known \cite{Cara2,DonTan1} that in this case $A$ is $s$-semigood for $s=d$.  In  contrast to this, when $A$ is not ``nearly square'', specifically, when $n>4\pi d$, $A$ can satisfy the condition $\bVSG_s(\xi,\theta,\rho,\sigma)$ only for $s\leq2$, no matter how large $\theta,\sigma,\rho$ are and how close to 1 $\xi<1$ is, see  {\color{blue} Online Supplement \ref{LemmaMoved}}.

\section{Verifiable sufficient conditions for $s$-semigoodness by Semidefinite Relaxation}
\label{sec:sdp}
Following d'Aspremont and El Ghaoui \cite{AspGaou}, we are about to derive another verifiable sufficient condition for $s$-semigoodness, now - via semidefinite relaxation.
The construction to follow is motivated by the development in the beginning of Section \ref{sec:ub}, according to which  $s$-semigoodness of $A$ is implied by the validity of (\ref{ub1}) for $\theta>1$ and  $\xi<1$.

Let, as before,
\[
\cX=\{x\in \Ker A: ~\Psi(x) \le 1\} ~\mbox{ and }~ \cU^{\xi,\theta}=\{V^{\xi,\theta}u:\;u\in \cU\},
\] where $\Psi$, $\cU$ and $V^{\xi,\theta}$ are defined in, respectively, \rf{psi2}, \rf{uu} and \rf{Vxite}.
The condition \rf{ub1} is equivalent to
\be
\max_{u,x}\left\{(V^{\xi,\theta}u)^Tx:\; u\in\cU,\; x\in\cX \right\}
\leq\xi.
\ee{sdp1}
Observe that for $u\in\cU$,  $x\in\cX$ the matrices $U=uu^T$,
$P=ux^T$ and $X=xx^T$ satisfy the relations
\begin{equation}\label{sdp2}
\begin{array}{ll}
\multicolumn{2}{l}{\exists t\in\bR^n,V\in\bS^{2n},\Lambda\in\bS^{2n}:}\\
(a)&\left[\begin{array}{c|c}U&P\cr\hline P^T&X\cr\end{array}\right]\succeq0;\\
(b)&\left\{\begin{array}{l}
U=\underbrace{\left[\begin{array}{c|c}I_n&-I_n\cr\end{array}\right]}_{:=L} \underbrace{\left[\begin{array}{c|c}V^{11}&V^{12}\cr \hline V^{12}&V^{11}\cr\end{array}\right]}_{:=V} L^T, \\0\leq V_{ij}\leq {1\over 2}, ~~V\succeq0, ~~V^{12}=[V^{12}]^T,\,\Tr(V)\leq s,\\
\sum_{i,j}V_{ij}\leq s^2, ~~V^{12}_{ij}=0~\forall i,j \in P_+;\\
\end{array}\right.\\
(c)&X=\underbrace{\left[\begin{array}{c|c|c|c}-I_p&0&\frac{1}{\theta}I_p&0\cr\hline 0&-I_q&0&I_q\cr\end{array}\right]}_{:=F}\Lambda F^T, ~~0\leq\Lambda_{ij}, ~~\Lambda\succeq0, ~~\sum_{i,j}\Lambda_{ij}\leq1;\\
(d_1)&\sum_{j\in P_+}\max[-P_{ij}, \theta P_{ij}] + \sum_{j\in P_n}|P_{ij}|\leq t_i, ~\forall i\in P_+,\\
(d_2)&\sum_j|P_{ij}|\leq t_i,~\forall i\in P_n,\\
(d_3)&t_i\leq 1~\forall i,~~\sum_i t_i\leq s;\\
(e)&AXA^T=0.\\
\end{array}
\end{equation}
Besides this,
$$
u^T(V^{\xi,\theta})^Tx=\Tr(V^{\xi,\theta}P^T).
$$
\begin{quote}Indeed, the latter relation, same as (\ref{sdp2}.$a$) and (\ref{sdp2}.$e$), is evident.  To verify (\ref{sdp2}.$b$), let $u_+=\max[u,0]$, $u_-=\max[-u,0]$, where $\max$ is acting coordinate-wise. Then
\bse
U&=&L\left[\begin{array}{c|c}u_+u_+^T&u_+u_-^T\cr\hline u_-u_+^T&u_-u_-^T\cr\end{array}\right]L^T=
L\left[\begin{array}{c|c}u_-u_-^T&u_-u_+^T\cr\hline u_+u_-^T&u_+u_+^T\cr\end{array}\right]L^T\\
&=&L\underbrace{\left[\begin{array}{c|c}{1\over2}[u_+u_+^T+u_-u_-^T]&{1\over 2}[u_+u_-^T+u_-u_+^T]\cr\hline {1\over 2}[u_-u_+^T+u_+u_-^T]&{1\over 2}[u_-u_-^T+u_+u_+^T]\cr\end{array}\right]}_{V}L^T,
\ese
and the matrix $V$ we have just defined clearly satisfies
all requirements from (\ref{sdp2}.$b$). To verify (\ref{sdp2}.$c$), observe that the extreme points of the set $\cX^+=\{x:\Psi(x)\leq1\}\supset\cX$
are the vectors $\pm e_i$, $i> p$, and $-e_i,\theta^{-1}e_i$, $i\leq p$, so that $x=F\lambda$ with $\lambda\in\bR^{2n}_+$, $\sum_i\lambda_i\leq1$; setting $\Lambda=\lambda\lambda^T$, we satisfy (\ref{sdp2}.$c$). To satisfy (\ref{sdp2}.$d$), it suffices to set $t_i=|u_i|$ for all $i$ and to take into account that $\max[-P_{ij},\theta P_{ij}]\geq |P_{ij}|$ for all $i,j$ due to $\theta\geq1$, and that $u_i\geq0$ for $i\in P_+$.
\end{quote}
It follows that a sufficient condition for (\ref{sdp1}) is
\begin{equation}\label{sdp3}
\Opt:=\max\limits_{\hbox{\scriptsize$\begin{array}{l}
     X,U\in\bS^n, \,V,\Lambda\in\bS^{2n}, \\
     P\in\bR^{n\times n},\,t\in\bR^n
    \end{array}$}}\left\{\Tr(V^{\xi,\theta}P^T): \hbox{\ (\ref{sdp2}) is satisfied}\right\}\leq\xi.
\end{equation}
The optimization problem in (\ref{sdp3}) clearly reduces to a semidefinite maximization program ${\cal S}$; by weak duality, the optimal value in the semidefinite dual ${\cal D}$ to ${\cal S}$ is $\geq \Opt$. It follows that the efficiently verifiable condition
$$
\Opt({\cal D})\leq\xi
$$
is a sufficient condition for $s$-semigoodness of $A$. Note that the above construction depends on $\theta\geq1$ and $\xi<1$ as parameters.
\paragraph{Remark.} Consider the case of $P_+=\emptyset$, where $\cX=\{x\in\bR^n:\|x\|_1\leq1,Ax=0\}\supset\cZ=\{x\in\bR^n:\|x\|_1\leq1\}$. In this case, the standard semidefinite relaxation of the set $\cC_*=\Conv\{xx^T:x\in\cZ\}$ is
\[
\cC=\left\{X:~X\succeq0,\sum_{i,j}|X_{ij}|\leq1\right\}
\](cf. \cite{AspGaou}). Note that (\ref{sdp2}.$c$) uses another semidefinite relaxation of
$\cC_*$, namely,
\[
\cC'=\left\{X:\exists \Lambda\in\bS^{2n}:\;
\begin{array}{l}\Lambda\succeq0,\Lambda_{i,j}\geq0\,\;\forall i,j,\;\sum_{i,j}\Lambda_{ij}\leq1\\
X=[I_n,-I_n]\Lambda[I_n,-I_n]^T\\ \end{array}\right\}.
\]
It is immediately seen that $\cC_*\subset \cC'\subset\cC$; a surprising fact is that the second of these inclusions is strict. Thus, the relaxation of $\cC_*$ given by $\cC'$ is less conservative than the standard relaxation given by $\cC$. As observed by A. d'Aspremont (private communication), the relaxation $\cC'$ can be further improved, namely, by replacing $ \cC'$ with
\[
\cC^+=\left\{X:\exists \Lambda=\left[\begin{array}{cc}\Lambda^{11}&\Lambda^{12}\cr\Lambda^{21}&\Lambda^{22}\cr\end{array}\right]
\in\bS^{2n}:\;
\begin{array}{l} \Lambda^{\mu\nu}\in\bR^{n\times n},\;\Lambda\succeq0,\;\Lambda_{i,j}\geq0\;\forall i,j\\
\;\sum_{i,j}\Lambda_{ij}\leq1, \; \Lambda^{12}_{ii}=0,\;1\leq i\leq n\\
X=[I_n,-I_n]\Lambda[I_n,-I_n]^T\\
\end{array}\right\}.\]
\par\noindent
{
Note that this idea can be used to improve the semidefinite relaxation given by $\cC$ as well. Specifically, the matrix $V$ as built in the justification of (\ref{sdp2}) clearly satisfies $(V^{12})_{ii}=0$, $1\leq i\leq n$, and we can add these linear constraints on $V$ to (\ref{sdp2}.$b$). Similarly, when representing a vector $x\in\cX^+$ as $F\lambda$ with $\lambda\in\bR^{2n}_+$, $\sum_i\lambda_i\leq1$, see the justification of (\ref{sdp2}), we clearly can ensure that $\lambda_i\lambda_{n+i}=0$, $1\leq i\leq n$, that is, the matrix $\Lambda$ we have built
in fact satisfies $\Lambda_{i,n+i}=\Lambda_{n+i,i}=0$, $1\leq i\leq n$, and we can add these linear constraints on $\Lambda$ to (\ref{sdp2}.$c$).

\section{Numerical results}
\label{sec:numerical}
In order to compare the performance of the proposed bounds on the maximal $s=s_*(A)$ for which a given matrix, $A$, is $s$-semigood, with the  bounds known from the literature, we present some preliminary numerical results for relatively small sensing matrices. Our goal is to see if the sign information on a signal allows to improve the bounds for $s_*(A)$ as compared to the bounds on the largest $s=s_0(A)$ for which $A$ is $s$-good.

We generate four sets of random matrices, which are normalizations (all columns scaled to be of $\|\cdot\|_2$-norm 1) of (a) Rademacher matrices (i.i.d. entries taking values $\pm1$ with probabilities 0.5), (b) Gaussian matrices (iid $\cN(0,1)$ entries), (c) Fourier matrices ---  $m\times n$ submatrices of the matrix of $n\times n$ Discrete Fourier Transform, and (d) Hadamard matrices --- $m\times n$ submatrices of the $n\times n$ Hadamard matrix\footnote{The Hadamard matrix $H_d$, $d=0,1,2,...$, has order $2^d\times 2^d$ and is given by the recurrence $H_0=1$, $H_{d+1}=[H_d,H_d;H_d,-H_d]$.}; in the cases (c,d), the $m$ rows comprising the submatrix were drawn at random from the $n$ rows of the ``parent'' matrix.  For each type, we set the number of columns to $n=256$ and vary the number of rows, $m=0.5n,\ldots,0.95n$.

We bound from below the value $s_0(A)$ using the  bound $s[\mu]$ by mutual incoherence and the bounds $s[\alpha_1]$ and $s[\alpha_s]$, computed through the LP-based verifiable sufficient conditions for $s$-goodness (see \cite[Section 6]{JNCS}).

{\em The lower bound on $s_*(A)$} is computed by invoking condition $\bVSG_s(\xi,\theta,\rho,\sigma)$, where  $\rho=\sigma=\infty$ and $\theta$ is set to once for ever fixed ``large enough'' value, and $\xi$ is set to 0.9999, see section \ref{LowerBounds} and Propositions \ref{propmain}, \ref{propnew}. Note that given a matrix $Y$, and setting $v=0$, one can compute the largest $s$ satisfying \eqref{condition} and thus ensuring the validity of $\bVSG_s(\xi,\theta,\rho,\sigma)$. We first compute the best lower bound $\underline{s}$ on $s_*(A)$ given by the $Y$-matrices generated when bounding $s_0(A)$. Then we compute the ``improved'' lower bound for $s_*(A)$ as follows: we check whether the condition $\bVSG_s(\xi,\theta,\rho,\sigma)$ holds true for $s=\underline{s}+1$, if it is the case, check whether this condition holds true for $s=\underline{s}+2$, and so on.

While the outlined lower bounds on $s_*(A)$ and $s_0(A)$ are efficiently computable via LP (when $\sigma=\rho=\infty$, the sufficient condition is easily checked by solving a Linear Programming program), the sizes of the resulting LPs are rather large. For instance, when $A$ is $m\times n$, the LP associated with \eqref{condition} has a $(2n^2+2n+1) \times ((m+2n)(n+1)+2)$ constraint matrix (compared to $(2n^2 + n) \times (n(m + n + 1) + 1)$ constraint matrices arising when computing lower bounds for $s_0(A)$). For instance, for $m = 230$ and $n = 256$, bounding $s_*(A)$ results in an LP program of the size $131,585\times190,696$, while computing a lower bound on $s_0(A)$ requires solving an LP problem of size $131,328\times124,673$. In all the computations, we used the state-of-the-art commercial LP solver {\tt mosekopt} \cite{Mosek}.

{\em The upper bounds} on $s_*(A)$ and on $s_0(A)$ are computed by the techniques from Section \ref{sec:ub} and \cite[Section 4.1]{JNCS}.

The results of our experiments and related CPU times  are presented in Table \ref{table:ExprResults}. The computations were carried out on a single core of an 8-core Intel Xeon E5520@2.27GHz CPU Linux workstation.

\begin{table}
\caption{Comparison of efficiently computable bounds on $s_*(A)$, $n=256$}
\label{table:ExprResults}
\centering
\scriptsize
{
\begin{tabular}{||c||c|cc|c|c|c||cc|c|c|c||}
\multicolumn{12}{c}{} \\
\multicolumn{12}{c}{\textbf{Fourier matrices}} \\
\hline
&  \multicolumn{4}{|c|}{\textbf{Unsigned}} & \multicolumn{2}{|c||}{\textbf{Nonnegative}} & \multicolumn{5}{|c||}{\textbf{CPU time (s)}}  \\
\hline
&  \multicolumn{3}{|c|}{LBs on $s_0(A)$} &  UB & LB & UB &
\multicolumn{3}{|c|}{Unsigned} & \multicolumn{2}{|c||}{Nonnegative} \\
\hline
$m$ & $s[\mu]$  &   $s[\alpha_1]$ & $s[\alpha_s]$ & $\bar{s}$ & $s[\alpha_s]$ & $\bar{s}$ & $s[\alpha_1]$ & $s[\alpha_s]$ & $\bar{s}$ & $s[\alpha_s]$ & $\bar{s}$  \\
\hline
\hline
128    &   3   &   5   &   5   &   12  &   5   &   47  &   0.8 &   1054.0  &   146.0    &   3114.4 &   172.9    \\
128 &   3   &   5   &   5   &   11  &   5   &   32  &   0.9 &   986.0   &   169.4    &   2891.5 &   311.5   \\
\hline
152 &   2   &   6   &   6   &   11  &   6   &   49  &   1.1 &   898.5   &   252.5    &   3680.2 &   179.6   \\
152    &   3   &   6   &   6   &   11  &   6   &   53  &   1.3 &   899.3   &   161.7    &   3836.7 &   183.5    \\
\hline
178    &   2   &   6   &   6   &   12  &   6   &   47  &   1.1 &   866.5   &   228.6    &   3976.0 &   294.0    \\
178 &   3   &   7   &   7   &   16  &   7   &   42  &   0.7 &   484.8   &   365.2    &   3216.8 &   416.9   \\
\hline
204 &   4   &   8   &   8   &   17  &   8   &   67  &   1.0 &   828.5   &   235.4    &   3829.7 &   209.2   \\
204    &   3   &   7   &   7   &   15  &   7   &   65  &   1.1 &   906.8   &   220.2    &   3914.4 &   197.4    \\
\hline
230 &   4   &   10  &   10  &   21  &   10  &   70  &   1.1 &   1879.9  &   300.5    &   4287.6 &   384.6   \\
230    &   4   &   9   &   9   &   20  &   9   &   65  &   1.0 &   856.6   &   286.5    &   4040.2 &   362.0    \\
\hline
242 &   5   &   11  &   11  &   26  &   11  &   89  &   1.7 &   1425.1  &   290.5    &   6444.1 &   513.0   \\
242    &   4   &   10  &   10  &   19  &   10  &   75  &   1.2 &   1920.6  &   265.3    &   4069.1 &   232.8    \\
\hline
\multicolumn{12}{c}{} \\
\multicolumn{12}{c}{\textbf{Hadamard matrices}} \\
\hline
&  \multicolumn{4}{|c|}{\textbf{Unsigned}} & \multicolumn{2}{|c||}{\textbf{Nonnegative}} & \multicolumn{5}{|c||}{\textbf{CPU time (s)}}  \\
\hline
&  \multicolumn{3}{|c|}{LBs on $s_0(A)$} &  UB & LB & UB &
\multicolumn{3}{|c|}{Unsigned} & \multicolumn{2}{|c||}{Nonnegative} \\
\hline
$m$ & $s[\mu]$  &   $s[\alpha_1]$ & $s[\alpha_s]$ & $\bar{s}$ & $s[\alpha_s]$ & $\bar{s}$ & $s[\alpha_1]$ & $s[\alpha_s]$ & $\bar{s}$ & $s[\alpha_s]$ & $\bar{s}$  \\
\hline
\hline
128    &   3   &   5   &   5   &   7   &   5   &   8   &   0.2 &   1148.1  &   77.8     &   3007.0 &   68.5     \\
128 &   2   &   5   &   5   &   7   &   5   &   7   &   0.3 &   1297.1  &   73.4     &   2894.4 &   116.8   \\
\hline
152 &   3   &   7   &   7   &   7   &   7   &   58  &   0.3 &   1224.4  &   47.9     &   3997.0 &   186.8   \\
152    &   4   &   7   &   7   &   13  &   7   &   58  &   0.2 &   1205.8  &   245.0    &   3962.6 &   310.4    \\
\hline
178    &   4   &   9   &   9   &   15  &   9   &   70  &   0.2 &   1269.8  &   238.9    &   4828.2 &   212.0    \\
178 &   4   &   9   &   9   &   15  &   9   &   19  &   0.3 &   1340.7  &   271.1    &   4923.3 &   342.8   \\
\hline
204 &   4   &   12  &   12  &   15  &   12  &   16  &   0.5 &   2908.1  &   131.2    &   6409.9 &   385.4   \\
204    &   5   &   12  &   12  &   15  &   12  &   16  &   0.4 &   2996.7  &   148.9    &   5507.9 &   253.9    \\
\hline
230 &   8   &   18  &   18  &   31  &   19  &   31  &   0.3 &   1860.1  &   250.8    &   9046.7 &   331.1   \\
230    &   8   &   18  &   18  &   31  &   18  &   39  &   0.4 &   2100.2  &   282.8    &   4081.3 &   396.8    \\
\hline
242 &   12  &   26  &   26  &   31  &   27  &   31  &   0.3 &   2015.1  &   92.7     &   7478.2 &   176.2   \\
242    &   12  &   26  &   26  &   31  &   26  &   31  &   0.3 &   1976.7  &   116.8    &   3597.9 &   412.0    \\
\hline
\multicolumn{12}{c}{} \\
\multicolumn{12}{c}{\textbf{Rademacher matrices}} \\
\hline
&  \multicolumn{4}{|c|}{\textbf{Unsigned}} & \multicolumn{2}{|c||}{\textbf{Nonnegative}} & \multicolumn{5}{|c||}{\textbf{CPU time (s)}}  \\
\hline
&  \multicolumn{3}{|c|}{LBs on $s_0(A)$} &  UB & LB & UB &
\multicolumn{3}{|c|}{Unsigned} & \multicolumn{2}{|c||}{Nonnegative} \\
\hline
$m$ & $s[\mu]$  &   $s[\alpha_1]$ & $s[\alpha_s]$ & $\bar{s}$ & $s[\alpha_s]$ & $\bar{s}$ & $s[\alpha_1]$ & $s[\alpha_s]$ & $\bar{s}$ & $s[\alpha_s]$ & $\bar{s}$  \\
\hline
\hline
128    &   1   &   5   &   5   &   14  &   5   &   53  &   27.8    &   1253.1  &    171.6   &  3388.7  &    124.8  \\
128 &   1   &   5   &   5   &   15  &   5   &   48  &   27.8    &   1361.5  &    191.1   &  3291.6  &   123.4    \\
\hline
152 &   2   &   6   &   6   &   18  &   7   &   65  &   38.4    &   1426.3  &    322.7   &  9592.1  &   136.3    \\
152    &   1   &   6   &   6   &   19  &   7   &   66  &   38.3    &   1183.0  &    218.9   &  9146.3  &    139.0  \\
\hline
178    &   2   &   7   &   8   &   25  &   9   &   78  &   44.2    &   2819.1  &    258.9   &  8032.1  &    225.8  \\
178 &   2   &   7   &   8   &   24  &   9   &   78  &   41.8    &   2481.7  &    256.0   &  8306.3  &   168.2    \\
\hline
204    &   2   &   10  &   11  &   32  &   12  &   92  &   51.1    &   1434.2  &    291.8   &  9738.5  &    209.3  \\
204 &   2   &   10  &   11  &   30  &   12  &   90  &   50.8    &   1316.6  &    448.3   &  9146.8  &   345.4    \\
\hline
230    &   2   &   14  &   16  &   41  &   19  &   107 &   61.8    &   2422.9  &    302.7   &  15235.2 &    162.2  \\
230 &   2   &   14  &   16  &   39  &   19  &   107 &   61.7    &   2466.2  &    624.0   &  15578.4 &   161.9    \\
\hline
242    &   2   &   20  &   23  &   47  &   27  &   116 &   64.8    &   3929.4  &    269.2   &  19828.7 &    178.1  \\
242 &   2   &   19  &   23  &   47  &   27  &   111 &   68.0    &   4242.4  &    277.8   &  20506.7 &   270.5    \\
\hline
\multicolumn{12}{c}{} \\
\multicolumn{12}{c}{\textbf{Gaussian matrices}} \\
\hline
&  \multicolumn{4}{|c|}{\textbf{Unsigned}} & \multicolumn{2}{|c||}{\textbf{Nonnegative}} & \multicolumn{5}{|c||}{\textbf{CPU time (s)}}  \\
\hline
&  \multicolumn{3}{|c|}{LBs on $s_0(A)$} &  UB & LB & UB &
\multicolumn{3}{|c|}{Unsigned} & \multicolumn{2}{|c||}{Nonnegative} \\
\hline
$m$ & $s[\mu]$  &   $s[\alpha_1]$ & $s[\alpha_s]$ & $\bar{s}$ & $s[\alpha_s]$ & $\bar{s}$ & $s[\alpha_1]$ & $s[\alpha_s]$ & $\bar{s}$ & $s[\alpha_s]$ & $\bar{s}$  \\
\hline
\hline
128 &   1   &   5   &   5   &   14  &   5   &   44  &   28.2    &   852.1   &    172.4   &  3283.2  &   114.7    \\
128    &   1   &   4   &   5   &   15  &   5   &   52  &   27.7    &   1913.9  &    177.7   &  3712.0  &    124.6  \\
\hline
152    &   2   &   6   &   6   &   19  &   7   &   58  &   35.4    &   981.0   &    214.1   &  8433.5  &    392.8  \\
152 &   1   &   6   &   6   &   19  &   7   &   58  &   38.9    &   1004.0  &    242.6   &  8231.7  &   373.3    \\
\hline
178 &   2   &   7   &   8   &   24  &   9   &   79  &   43.0    &   2164.4  &    393.9   &  10294.7 &   368.2    \\
178    &   2   &   7   &   8   &   25  &   9   &   77  &   47.6    &   2390.3  &    263.1   &  9548.8  &    374.0  \\
\hline
204 &   2   &   10  &   11  &   32  &   12  &   88  &   58.0    &   1363.6  &    293.3   &  11496.7 &   274.1    \\
204    &   2   &   10  &   11  &   32  &   12  &   91  &   51.7    &   1218.4  &    293.4   &  12497.2 &    529.5  \\
\hline
230 &   2   &   14  &   17  &   41  &   19  &   102 &   70.4    &   3200.9  &    339.7   &  18771.3 &   431.6    \\
230    &   2   &   14  &   16  &   39  &   19  &   106 &   61.5    &   2118.4  &    485.4   &  18959.5 &    435.0  \\
\hline
242    &   2   &   19  &   22  &   46  &   27  &   113 &   73.6    &   2212.8  &    277.4   &  26874.6 &    269.2  \\
242 &   2   &   20  &   23  &   47  &   27  &   112 &   65.3    &   2995.2  &    426.7   &  21308.7 &   191.7    \\
\hline
\end{tabular}
}
\end{table}

The results in Table \ref{table:ExprResults} merit some comments. We observe that our LP-based efficiently computable lower bounds  on $s_0(A)$ and $s_*(A)$ clearly outperform the bounds based on mutual incoherence.
We notice that for Fourier and Hadamard matrices, the lower bounds on $s_*(A)$ and $s_0(A)$ are nearly always the same, except for two Hadamard instances with $m=230$ and $m=242$. On the other hand, for Gaussian and Rademacher matrices, as the number of rows $m$ approaches the number of columns $n$, the difference between the best certified lower bounds on $s_*(A)$ and on $s_0(A)$ increases (for the sizes we have considered, this difference attains $5$ for the Gaussian matrix with $m=242$).
While for Gaussian, Rademacher and Fourier matrices, the upper bounds on $s_*(A)$ become loose (they are twice or three times higher than the upper bounds on $s_0(A)$), these bounds become tighter in the case of Hadamard matrices. Further, for some matrices  the lower and the upper bound on $s_0(A)$ match (e.g., the Hadamard matrix with $m=152$), what allows to identify the exact value of $s_0(A)$ . Moreover, we have observed samples of smaller random Hadamard matrices (with $n=128$) for which the lower bounds and upper bounds on both $s_*(A)$ and $s_0(A)$ coincide, which implies $s_*(A)=s_0(A)$  in these cases.

\section{Matching pursuit algorithm}
\label{sec:mp}
The Matching Pursuit algorithm for signal recovery has been first introduced in \cite{malzhang} and is motivated by the desire to provide a reduced complexity alternative to the $\ell_1$-recovery problem. Several implementations of Matching Pursuit has been proposed in the Compressive Sensing literature (see, e.g., the review \cite{DonEl}). All of them are based on successive Euclidean projections of the signal and the corresponding performance results rely upon the bounds on mutual incoherence $\mu(A)$ of the sensing matrix. We are about to show that the LP-based verifiable sufficient conditions from the previous section can be used to construct a specific version of the Matching Pursuit algorithm which we refer to as {\em Non-Euclidean Matching Pursuit (NEMP) algorithm}.

Suppose that we have in our disposal $\tau, \tau_{\pm} \geq 0$ and a matrix $Y=[y_1, ... ,y_n]$, such that
\begin{equation}\label{feq99}
\begin{array}{lll}
(a)& -\tau_-\leq  [I-Y^TA]_{ij}\leq \tau_+,&\forall i\in P_+, ~\forall j,\\
(b)& -\tau \leq [I-Y^TA]_{ij}\leq \tau,&\forall i\in P_n,~\forall j,\\
(c)&\|y_j\|_*\leq\sigma,&\forall j.
\end{array}
\end{equation}

Consider a signal {$w\in\cP_n$} such that $\|w-w^s\|_1\le \mu$, where $w^s$ is the vector obtained from $w$ by replacing all but $s$ largest magnitudes of entries in $w$ with zeros, and let $y$ and $\delta$ be such that $\|Aw-y\|\leq\delta$.
\par
Suppose that
\begin{equation}\label{rho}
\rho=s\max\{\tau_+,\tau_-,\tau\}<1.
\end{equation}
To simplify notation, we denote $\max[a,b]$ by $a \vee b$. Consider the following iterative procedure:
\begin{algorithm}\label{MP}
$~$
\begin{enumerate}
\item \underline{Initialization:} Set $v^{(0)}=0$, $
    \alpha_0={\|Y^Ty\|_{s,1}+s\sigma\delta+\mu\over 1-\rho}.$
\item
\underline{Step $k$, $k=1,2,...$:} Given $v^{(k-1)}\in \bR^n$ and $\alpha_{k-1}\geq 0$, compute
\begin{enumerate}
\item $u=Y^T(y-Av^{(k-1)})$ and $n$ segments
\[
S_i=\left\{\begin{array}{ll}
~[u_i-\tau_-\alpha_{k-1}-\sigma\delta,\,u_i+\tau_+\alpha_{k-1}+\sigma\delta], & i\in P_+,\\
~[u_i-\tau\alpha_{k-1}-\sigma\delta,\,u_i+\tau\alpha_{k-1}+\sigma\delta], & i \in P_n.
\end{array}
\right.
\]
Define $\Delta\in\bR^n$ by setting
\bse
\Delta_i=
\left\{
\begin{array}{ll}
~~[u_i-\tau_- \alpha_{k-1}-\sigma\delta]_+, & i\in P_+,\\
~~[u_i-\tau \alpha_{k-1}-\sigma\delta]_+,& i\in P_n,\;\;u_i\ge 0,\\
-[|u_i|-\tau \alpha_{k-1}-\sigma\delta]_+, & i\in P_n,\;\;u_i<0
\end{array}\right.
\ese
(here $[a]_+=\max[0,a]$).
\item
Set  $v^{(k)} = v^{(k-1)} + \Delta$ and
\begin{equation}\label{finitedif}
\alpha_k=s[ 2\tau\vee (\tau_-+\tau_+)]\alpha_{k-1}+2s\sigma\delta+\mu.
\end{equation}
and loop to step $k+1$.
\end{enumerate}
\item The approximate solution found after $k$ iterations is $v^{(k)}$.
\end{enumerate}
\end{algorithm}

\begin{proposition}\label{Prop_NEMP} Assume that $w_i\geq0$ for $i\in P_+$, {\rm (\ref{rho})} takes place, and that $\|w-w^s\|_1\leq\mu$ with a known in advance value of $\mu$. Then
the approximate solution $v^{(k)}$ and the value $\alpha_k$ after the $k$-th step of   Algorithm \ref{MP} satisfy
$$
\begin{array}{llcll}
(a_{k})&\mbox{for all $i$}\;\;v^{(k)}_i\in\Conv\{0;w_i\},&&
(b_{k})&\|w-v^{(k)}\|_1%=\sum_{i\in P_+}[w_i-v^{(k)}_i]+\sum_{i\in P_n}|w_i-v^{(k)}_i|
\leq\alpha_{k}.
\end{array}
$$
\end{proposition}
For proof, see Appendix \ref{App:Prop_NEMP}.

Let
$$
\lambda=s[ 2\tau\vee (\tau_-+\tau_+)];
$$
if $\lambda<1$, then also $\rho<1$, so that Proposition \ref{Prop_NEMP} holds true. Furthermore, by (\ref{finitedif}) the sequence $\alpha_k$ converges exponentially fast to the limit $\alpha_\infty:={2s\sigma\delta+\mu\over 1-\lambda}$:
\bse
\alpha_k=\lambda^k[\alpha_0-\alpha_\infty] +\alpha_\infty.
\ese
Note that when $P_+=\emptyset$, we can set $\tau_-=\tau_+=0$ to obtain $\lambda=2s\tau$; in the case of $P_n=\emptyset$, by setting $\tau=0$, we have $\lambda=s(\tau_-+\tau_+)$.

The bottom line is: if the optimal value in the convex program
$$
\Opt=\min\limits_{\tau,\tau_\pm,Y}\left\{ s[ 2\tau\vee (\tau_-+\tau_+)]:\;\begin{array}{cl}-\tau_-\leq  [I-Y^TA]_{ij}\leq \tau_+,&\forall i\in P_+, ~\forall j\\
-\tau\leq  [I-Y^TA]_{ij}\leq \tau,&\forall i\in P_n, ~\forall j\\
\tau,\tau_\pm\geq0\\
\end{array}\right\}
$$
is $<1$, the above procedure, as yielded by an optimal solution to the latter problem,  possesses the following properties:
\begin{enumerate}
\item All approximations $v^{(k)}$, $k=0,1,...$ of $w$ are supported on the support of $w$;
\item For $i\in P_+$, $v^{(k)}_i\geq0$ are nondecreasing in $k$ and are $\leq w_i$ for all $k$;
\item For $i\in P_n$,
\begin{itemize}
\item if $w_i>0$, then $0\leq v_i^{(k)}\leq w_i$ and $v_i^{(k)}$ are nondecreasing in $k$;
\item if $w_i<0$, then $w_i\leq v_i^{(k)}\leq 0$ and $v_i^{(k)}$ are nonincreasing in $k$;
\end{itemize}
\item As $k$ grows, the upper bound $\alpha_k$ on the $\ell_1$-error of approximating $w$ by $v^{(k)}$ goes exponentially fast
to \[
\alpha_\infty={2s\sigma\delta+\mu\over 1-\Opt}.
\]
\end{enumerate}

Let now $\xi\in[0,1)$, $\sigma\ge 0$ and $\theta\ge1$ and suppose that an $m\times n$ matrix $A$ satisfies the following condition:
\begin{quote}
$\overline{\bVSG}_s(\xi,\sigma,\theta)$: {\sl There exists $m\times n$ matrix $Y=[y_1,...,y_n]$ such that $\|y_i\|_*\leq\sigma$ for all $i$ and}
\begin{equation}\label{ovlA}
\begin{array}{ll}
-{\xi\over {(1+\xi)s}}\leq [I-Y^TA]_{ij}\leq{\xi\over {(1+\xi)s}} &\forall i \not \in P_+, ~\forall j, \\
-{\xi\over {(1+\xi\theta)s}}\leq [I-Y^TA]_{ij}\leq{\xi\over {(1+\xi\theta)s}} &\forall i \in P_+, ~\forall j \not \in P_+, \\
-{\xi\over {(1+\xi\theta)s}}\leq [I-Y^TA]_{ij}\leq{\xi\theta\over {(1+\xi\theta)s}} &\forall i,j \in P_+.
\end{array}
\end{equation}
\end{quote}
Observe that \rf{ovlA} is a system of convex inequalities in $Y$. Further,
$\overline{\bVSG}_s(\xi,\sigma,\theta)$ certainly implies $\bVSG_s(\xi,\theta,0,\sigma)$, and is therefore sufficient condition for $s$-semigoodness of the matrix $A$.
\par
When $\overline{\bVSG}_s(\xi,\sigma,\theta)$ is satisfied with $\xi\in(0,1)$ and $\theta>1$, by taking
\[
\tau_-={\xi\over {(1+\xi\theta)s}}, \;\;\;\tau_+={\xi\theta\over {(1+\xi\theta)s}}\;\;\mbox{ and}\;\;\tau={\xi\over {(1+\xi)s}},
\]
we obtain
\begin{equation}\label{lambda}
\lambda=\max\left({\xi+\xi\theta\over 1+\xi\theta}, \;{2\xi\over 1+\xi} \right)<1.
\end{equation}
Combining this condition with Proposition \ref{Prop_NEMP} gives:

\begin{corollary}\label{MPB} Suppose that $A$ satisfies the condition $\overline{\bVSG}_{s}(\xi,\sigma,\theta)$ with certain $\xi\in(0,1)$, $\sigma\ge 0$ and $\theta\ge1$. Let {$w\in\cP_n$} be a  vector with $\|w-w^s\|_1\le \mu$ where $w^s$ is the vector obtained from $w$ by replacing all but $s$ largest in magnitude entries in $w$ with zeros, and let $y$ be such that $\|Aw-y\|\leq \delta$.
Then the approximate solution $v^{(t)}$ found by  Algorithm \ref{MP} after $t$ iterations satisfies $v^{(t)}_i \geq 0$ for all $i \in P_+$ and
\[
\|w-v^{(t)}\|_1\leq {2s\sigma\delta+\mu\over1-\lambda}+\lambda^{t}\left[{\|Y^Ty\|_{s,1}+s\sigma\delta+\mu\over1-\rho}-{2s\sigma\delta+\mu\over 1-\lambda}\right],
\]
where $\lambda$ is given by {\rm (\ref{lambda})} and $\rho={\xi\theta\over 1+\xi\theta}$.
\end{corollary}

It should be noted the NEMP algorithm has several drawbacks as compared with the $\ell_1$-recovery. First, the pursuit algorithm requires a priori knowledge of several parameters ($\sigma$, $Y$, $\tau$, $\tau_-$, $\tau_+$, $s$ and $\mu$). Second,  the value $(1-\lambda)^{-1}(2s\sigma\delta+\mu)$ is a conservative upper bound on the error of the $\ell_1$-recovery, but the error bound in Corollary \ref{MPB} is exact.
On the other hand, the NEMP algorithm can be an interesting option if the $\ell_1$-recovery is to be used repeatedly on the observations obtained with the same sensing matrix $A$; the numerical complexity of the pursuit algorithm for a given matrix $A$ may only be a fraction of that of the $\ell_1$-recovery, especially when used on high-dimensional data.
\par
Our concluding remark is on the condition
\begin{equation}\label{mui}
{\mu(A)\over 1+\mu(A)}<{1\over 2s},
\end{equation} where $\mu(A)$ is the mutual incoherence of $A$ (see (\ref{muin})). This condition is usually used in order to establish convergence results for the Matching Pursuit algorithms  (see, e.g. \cite{DonElTem,Elad,ElZib}). As it is immediately seen, when $\mu(A)$ is well defined (i.e., all columns in $A$ are nonzero), the matrix $Y=[y_1,...,y_n]$ with the columns
\[
y_i={A_i\over (1+\mu(A))A_i^TA_i}
 \]satisfies for all $i=1,...,m$ and $j=1,...,n$ the relations
 \[
 |[I-Y^TA]_{ij}|\leq {\mu(A)\over 1+\mu(A)}.
 \]
 In the case of (\ref{mui}), setting $\theta=1$ and specifying $\xi$ from the relation ${\xi\over 1+\xi}={s\mu(A)\over 1+\mu(A)}$, we get $0<\xi<1$ and meet all inequalities in (\ref{ovlA}). It follows that $Y$ certifies the validity of the condition $\overline{\bVSG}_s(\xi,\sigma,1)$ with the outlined $\xi$ and with all $\sigma\geq\max\limits_i{\|A_i\|_*\over(1+\mu(A))\|A_i\|_2^2}$,
and thus the above $Y$ can be readily used in Matching Pursuit. Note that in the situation in question Corollary \ref{MPB} recovers some results from  \cite{DonElTem,Elad,ElZib}.

\appendix
\section{Proof of Proposition
\protect{\ref{prop1}}} \label{App:Prop1}

{\bf }
{\bf (i)$\Rightarrow$(ii):}
Let $A$ be $s$-semigood, and let, in contrast to what is stated by (ii), $J$ be a subset of $\{1,...,n\}$ with $\Card(J)\leq s$ and $x\in\Ker A\backslash \{0\}$ be such that $x_i\leq0$ for all $i\in P_+ \setminus J$ and \[\sum_{i\in J\cap P_+}x_i + \sum_{i\in J\cap P_n}|x_i| \geq \sum_{i\not\in J}|x_i|.
 \]Let $I=(J\cap P_n) \cup \{i \in J\cap P_+: x_i \ge 0\}$ so that $I\subseteq J$. From the construction of $I$, we have $x_i\leq0$ for $i\in J \setminus I$ implying that $x_i\leq0$ for $i\in P_+ \setminus I$. Further,
\begin{align*}
\sum_{i\in I\cap P_+}x_i + \sum_{i\in I\cap P_n}|x_i|
&= \sum_{i\in J\cap P_+}x_i - \sum_{i\in J\setminus I}x_i + \sum_{i\in J\cap P_n}|x_i| \\
&{\geq} \sum_{i\not\in J}|x_i| - \sum_{i\in J\setminus I}x_i
= \sum_{i\not\in J}|x_i| + \sum_{i\in J\setminus I}|x_i | = \sum_{i\not\in I}|x_i|.
\end{align*}
Hence $I$ also violates the condition in (ii). Setting $u_i=x_i$ when $i \in I$ and $u_i=0$ otherwise and setting $v=u-x$, we have $u_i\geq0$ for any $i \in I\cap P_+$, $u_i=0$ for any $i \in P_+ \setminus I$, and $v_i\geq0$ for $i \in P_+ \setminus I$, $v_i=0$ for $i \in I\cap P_+$ and $\sum_i|u_i|\geq\sum_i|v_i|$. In addition, $Au=Av$ due to $Ax=0$, and $u$ is $s$-sparse; finally, $u\neq v$ due to $x\neq0$. We see that the $s$-sparse vector $u\in\cP_n$ is not the unique solution to \[
\min_z\left\{\sum_i|z_i|: \;Az=Au, ~~z_i \geq 0 ~\forall i \in P_+\right\},\]
 which is a desired contradiction.
\par
{\bf (ii)$\Rightarrow$(iii):} Let $A$ satisfy (ii). Let ${\cal J}$ be the family of all subsets $J$ of $\{1,...,n\}$ of cardinality $\leq s$.  For $J\in{\cal J}$, let \[
X_J=\{x\in \Ker A: \|x\|_1=1, ~x_i\leq0\ ~\forall i\in P_+ \setminus J\}.\]
 Assuming that $X_J\neq\emptyset$, let $x\in X_J$. By (ii), we have
\[
%begin{equation}\label{neq10}
%\forall x\in X_J:
\sum_{i\in J\cap P_+}x_i + \sum_{i\in J\cap P_n}|x_i| < \sum_{i\not\in J}|x_i|.
\]
We claim that $\sum_{i\not\in J}|x_i|>0$.
 \begin{quote}
 Indeed, otherwise $x_i \neq 0$ implies that $i\in J$. Let $I_+$ and $I_-$ be the subsets of $J$ such that $x_i>0$ for $i\in I_-$ and $x_i<0$ for $i\in I_+$. At least one of these sets is nonempty due to $x\neq0$. W.l.o.g. we can assume that $\sum_{i\in I_+}x_i\geq\sum_{i\in I_-}|x_i|$ (otherwise we could replace $x$ with $-x$ and swap $I_+$ and $I_-$). Applying (ii) to $x$ and to $I_+$ in the role of $J$, we should have
 $$
 \sum_{i\in I_+\cap P_+}x_i + \sum_{i\in I_+\cap P_n}|x_i| = \sum_{i\in I_+}x_i<\sum_{i\not\in I_+}|x_i|=\sum_{i\in I_-}|x_i|,
 $$
 which is not the case. This contradiction shows that $\sum_{i\not\in J}|x_i|>0$ whenever $x\in X_J$.
 \end{quote}

From our claim it follows that the function
$$
\sum_{i\in J\cap P_+}x_i + \sum_{i\in J\cap P_n}|x_i| \over \sum_{i\not\in J}|x_i|
$$
is continuous on $X_J$ and is $<1$ at every point of this set. Since $X_J$ is compact, we conclude that when $J\in{\cal J}$ is such that $X_J\neq\emptyset$, there exists $\xi_J<1$ such that
$$
\sum_{i\in J\cap P_+}x_i + \sum_{i\in J\cap P_n}|x_i| \leq \xi_J \sum_{i\not\in J}|x_i| \mbox{ for any }x\in X_J.
$$
Setting $\xi=\max\limits_{J\in{\cal J}:X_J\neq\emptyset}\xi_J$, we clearly ensure the validity of (iii). The implication (ii)$\Rightarrow$(iii) is proved.
\par
\par
{\bf (iii)$\Rightarrow$(i):} Let (iii) take place; let us prove that $A$ is $s$-semigood. Thus, let $u$ with $u_i\geq0$ for all $i \in P_+$ be $s$-sparse; we should prove that $u$ is the unique optimal solution to the problem \[
\min\limits_z\left\{\sum_i|z_i|:Az=Au, ~z_i\geq0~\forall i \in P_+\right\}.
\]
Assume, on the contrary to what should be proved, that the latter problem has an optimal solution $v$ different from $u$, and let $x=u-v$, so that $x\in\Ker A$ and $x\neq0$. Setting ${I=\{i:u_i\neq0\}}$, we have $\Card(I)\leq s$ and $x_i\leq0$ when $i\in P_+ \setminus I$,  whence by (iii)
$$
\sum_{i\in I\cap P_+}x_i + \sum_{i\in I\cap P_n}|x_i| \leq \xi \sum_{i\not\in I}|x_i|=\xi \sum_{i\not\in I}|v_i|,
$$
whence also
\begin{equation}\label{neq11}
\underbrace{\sum_{i\in I\cap P_+}u_i + \sum_{i\in I\cap P_n}|u_i|}_{=\sum_{i\in I}|u_i|}
\leq \underbrace{\sum_{i\in I\cap P_+}v_i + \sum_{i\in I\cap P_n}|v_i|}_{=\sum_{i\in I}|v_i|} + \xi\sum_{i\not\in I}|v_i|.
\end{equation}
Since $\sum_i|v_i| \leq \sum_i|u_i| = \sum_{i\in I}|u_i|$ due to the origin of $v$, (\ref{neq11}) implies that $\sum_{i\not\in I}|v_i|=0$,
that is, both $u$ and $v$ are supported on $I$, { so that $x$ is supported on $I$ as well. Now let $I_+=\{i\in I\cap P_+:x_i\geq0\}$,
$I_-=\{i\in I\cap P_+:x_i<0\}$ and $I_n=I\cap P_n$. Replacing, if necessary, $x$ with $-x$ and swapping $I_+$ and $I_-$, we can assume that
$\sum_{i\in I_+}x_i=\sum_{i\in I_+}|x_i|\geq \sum_{i\in I_-}|x_i|$. Applying (iii) to $x$ and to $I_+\cup I_n$ in the role of $J$, we get
$$
\sum_{i\in I_+}x_i+\sum_{i\in I_n}|x_i|\leq\xi\sum_{i\in I_-}|x_i|,
$$
thereby $\sum_{i\in I_+}x_i=\sum_{i\in I_n}|x_i|=\sum_{i\in I_-}|x_i|=0$ due to $\sum_{i\in I_+}x_i\geq\sum_{i\in I_-}|x_i|$. Thus, $x=0$, which is a desired contradiction.}
\par We have proved that the properties (i) -- (iii) of $A$ are equivalent to each other.
\par
{\bf (iii)$\Leftrightarrow$(iv):} The implication (iv)$\Rightarrow$(iii) is evident. Let us prove the inverse implication. Thus, let $A$ satisfy {(iii)} (and thus -- (i) -- (ii) as well), and let $\xi'\in(\xi,1)$. Let, as above, ${\cal J}$ be the family of all subsets $J$ of $\{1,...,n\}$ of cardinality $\leq s$. Let $X=\{x\in\Ker A: \|x\|_1=1\}$, and let $J\in {\cal J}$. Let $x\in X$. We claim that there exists a neighborhood $U_x$ of $x$ in $X$ and $\theta_{J,x}\in[1,\infty)$ such that for any $u\in U_x$ and $\theta\geq\theta_{J,x}$ it holds
\begin{equation}\label{neq7}
\sum_{i\in J\cap P_+}u_i + \sum_{i\in J\cap P_n}|u_i| \leq \xi'
\left(  \sum_{i\in P_n\setminus J} |u_i| + \sum_{i\in P_+ \setminus J}\max[-u_i,\theta u_i] \right).
\end{equation}
\begin{quote} The claim is clearly true when there exists ${ i\in P_+\setminus J}$ such that $x_i>0$. Now assume that $x_i\leq0$ for ${ i\in P_+\setminus J}$. Then $\sum_{i\not\in J}|x_i|>0$. Indeed, otherwise $x_i=0$ for all $i\not\in J$, which combines with $s$-semigoodness of $A$ and the relation $Ax=0$ to imply that $x=0$ (since assuming $x\neq0$, we have $x=u-v$ with $s$-sparse $u\geq0,v\geq0$ with non-overlapping supports, and $Au=Av$ due to $Ax=0$, which of course contradicts the $s$-semigoodness of $A$), while $x$ definitely is nonzero (since $\|x\|_1=1$ due to $x\in X$). Now, since $x\in\Ker A$ and $x_i\leq 0$, ${ i\in P_+\setminus J}$, we have
$$
\sum_{i\in J\cap P_+}x_i + \sum_{i\in J\cap P_n}|x_i| \leq \xi \sum_{i\not\in J} |x_i|<\xi'\sum_{i\not\in J}|x_i|
$$
where the first inequality is due to (iii), and the second -- due to $\sum_{i\not\in J}|x_i|>0$. The concluding strict inequality
clearly implies the validity of (\ref{neq7}) with $\theta=1$, provided that $U_x$ is a small enough neighborhood of $x$.
Thus, our claim is true.\end{quote}
From the validity of our claim, extracting from the covering $\{U_x\}_{x\in X}$ of the compact set $X$ a finite subcovering,
we conclude that there exists $\theta_J\in[1,\infty)$ such that
$$
\forall (x\in X,\;\theta\geq\theta_J): \sum_{i\in J\cap P_+}x_i + \sum_{i\in J\cap P_n}|x_i| \leq \xi' \left(  \sum_{i\in P_n\setminus J} |x_i| + \sum_{i\in P_+ \setminus J}\max[-x_i,\theta x_i] \right).
$$
Setting $\theta=\max_{J\in{\cal J}}\theta_J$, we see that $A$ satisfies $\bSG_s(\xi',\theta)$.
\par
{\bf (iv)$\Rightarrow$(v):} Let $A$ satisfy $\bSG_s(\xi,\theta)$ for certain $\xi\in(0,1)$, $\theta\in[1,\infty)$ and let $\|\cdot\|$ be a norm on $\bR^m$. Let, further, $P$  be the orthogonal projector of $\bR^n$ on $\Ker A$. Then clearly with a properly chosen $C$ one has
$$
\|Px-x\|_1\leq C\|Ax\|
$$
for any $x\in \bR^n$. Now let $J$ be a subset of $\{1,...,n\}$ of cardinality $\leq s$, $x\in \bR^n$ and $u=Px$. We have
\bse
\lefteqn{\sum_{i\in J\cap P_+}x_i + \sum_{i\in J\cap P_n}|x_i| \leq \sum_{i\in J\cap P_+}u_i + \sum_{i\in J\cap P_n}|u_i| +\sum_{i\in J}|u_i-x_i|}\\
&\leq& \xi \left[  \sum_{i\in P_n\setminus J} |u_i| + \sum_{i\in P_+ \setminus J}\max[-u_i,\theta u_i] \right] +\sum_{i\in J}|u_i-x_i|\\
&\leq&{\xi\left[\sum_{i\in P_n\setminus J}[|x_i|+|u_i-x_i|]+\sum_{i\in P_+\setminus J}[{\max}[-x_i,\theta x_i]+\theta|x_i-u_i|]\right]+\sum_{i\in J}|u_i-x_i|}\\
&\leq&
\xi\left[\sum_{i\in P_n\setminus J}|x_i|+\sum_{i\in P_+\setminus J}\max[-x_i,\theta x_i]\right]+\max[1,\theta\xi]\|x-u\|_1\\
&\leq& \xi\left[\sum_{i\in P_n\setminus J}|x_i|+\sum_{i\in P_+\setminus J}\max[-x_i,\theta x_i]\right]+\max[1,\theta\xi]C\|Ax\|,\\
\ese
so that $A$ satisfies $\bSG_{s,\beta}(\xi,\theta)$ with $\beta=\max(1,\theta\xi)C$. The implication (iv)$\Rightarrow$(v) is proved.
\par
{\bf (v)$\Rightarrow$(vi)$\Rightarrow$(iii):} These implications are evident. \qed

\section{Proof of Proposition \protect{\ref{IRB}}}\label{App:IRB}

{\bf }
Let $I$ be the support of $w^s$, $\bar{I}$ be the complement of $I$ in $\{1,...,n\}$, and let $z=w-x$. We denote $I_+=\{i\in I: z_i \ge 0\}$, $\bar{I}_+=\{i\in \bar{I}: z_i \ge 0\}$, and $I_-=I \setminus I_+$, $\bar{I}_-=\bar{I} \setminus \bar{I}_+$. Observe that $w$ is a feasible solution to (\ref{problem!}), so that
\begin{equation}\label{Opt!}
\|x\|_1\leq \|w\|_1+\nu.
\end{equation}
Obviously, $|x_i|-|w_i|\ge -|z_i|$ and $|x_i|-|w_i|\ge |z_i|-2|w_i|$. Now
using $x_i, w_i \ge 0 ~\forall i \in P_+$, and $z_i \ge 0 ~\forall i \in I_+$, we get
\bse
\nu &\ge& \sum_i[|x_i|-|w_i|]\qquad\hbox{[by (\ref{Opt!})]}\\
&\ge&\sum_{i \in I_+ \cap P_+} \underbrace{(x_i-w_i)}_{=-z_i} + \sum_{i \in I_- \cap P_+} \underbrace{(x_i-w_i)}_{=-z_i = |z_i|} + \sum_{i \in \bar{I}_- \cap P_+} \underbrace{(x_i-w_i)}_{=-z_i=|z_i|} +\sum_{i \in \bar{I}_+ \cap P_+} \underbrace{(x_i-w_i)}_{=-z_i\ge-w_i}\\
&&+ \sum_{i \in P_n} (|x_i|-|w_i|)\\
&\ge &-\sum_{i \in I_+ \cap P_+}z_i+\sum_{i \in I_- \cap P_+}|z_i|
+\sum_{i \in \bar{I}_- \cap P_+}|z_i|-\sum_{i \in \bar{I}_+ \cap P_+}w_i\\
&&-\sum_{i\in I\cap P_n}|z_i|+\sum_{i\in \bar{I}\cap P_n}(|z_i|-2|w_i|),
\ese
or, equivalently,
\begin{equation}\label{ai1}
\begin{array}{l}
\sum_{i \in I_- \cap P_+} |z_i|
+\sum_{i \in \bar{I}_-\cap P_+}|z_i|+\sum_{i \in \bar{I}\cap P_n} |z_i| \\
\le \nu + \sum_{i \in I_+ \cap P_+} z_i + \sum_{i \in I \cap P_n} |z_i|
+\sum_{i\in \bar{I}_+\cap P_+} w_i+2\sum_{i \in \bar{I}\cap P_n} |w_i|.
\end{array}
\end{equation}

On the other hand, we have
\begin{equation}\label{wehave22}
\|Az\|=\|Aw-Ax\|\le \|Aw-y\|+\|Ax-y\|\le \e+\delta.
\end{equation}
Then by condition $\bSG_{s,\beta}(\xi,\theta)$ with $(I_+ \cap P_+) \cup (I \cap P_n)$ in the role of $J$, we get
\begin{equation}\label{eqprev}
\begin{array}{l}
\underbrace{\sum_{i\in I_+ \cap P_+} z_i + \sum_{i\in I \cap P_n} |z_i|}_{:=\kappa}
\le \beta\|Az\| + \xi  \left[\sum_{i \in \bar{I} \cap P_n} |z_i| + \sum_{i \in (\bar{I}\cap P_+) \cup (I_- \cap P_+)}\psi(z_i) \right] \\
\kappa \le \beta\|Az\| + \xi \bigg[\underbrace{{{\sum}_{i\in \bar{I} \cap P_n}|z_i|}+ {\sum}_{i\in I_- \cap P_+}|z_i|+{{\sum}_{i\in \bar{I}_-\cap P_+}|z_i|} +\theta{\sum}_{i\in \bar{I}_+ \cap P_+}z_i}_{:=\tau(\theta)}\bigg]\\
\end{array}
\end{equation}
Let us derive a bound on $\tau(\theta)$. Now \rf{ai1} implies, independently of whether $\bSG_{s,\beta}(\xi,\theta)$ is or is not true, the first inequality in the following chain:
\be
\tau(\theta)&\le &
\nu + \sum_{i \in I_+ \cap P_+} z_i + \sum_{i \in I \cap P_n} |z_i|
+\sum_{i\in \bar{I}_+\cap P_+} w_i+2\sum_{i \in \bar{I}\cap P_n} |w_i|+
\theta\sum_{i\in \bar{I}_+\cap P_+}z_i\nn
&\le &\nu + \kappa
+(1+\theta)\sum_{i\in \bar{I}_+\cap P_+} w_i+2\sum_{i \in \bar{I}\cap P_n} |w_i|\qquad\hbox{\ [since $w_i\geq z_i$ for $i\in P_+$]}\nn
&\le &\nu+\kappa+(1+\theta)\mu, \qquad\hbox{\ [since $\theta\geq1$ and $\sum_{i\in \bar{I}}|w_i|\leq\mu$]},
\ee{ai2}
and, in particular,
\be
\tau(1) = \sum_{i\in I_- \cap P_+}|z_i|+\sum_{i\in \bar{I}}|z_i| \le \nu + \kappa + 2 \mu.
\ee{ai22}
Combining (\ref{wehave22}), (\ref{eqprev}) and \rf{ai2}, we obtain
$$
\kappa\leq\beta(\e+\delta)+\xi\left[\nu+\kappa+(1+\theta)\mu\right],
$$
and thereby,
\[
\kappa=\sum_{i\in I_+ \cap P_+} z_i + \sum_{i\in I \cap P_n} |z_i| \le \frac{\beta(\e+\delta)+\xi(\nu+(\theta+1)\mu)}{1-\xi} .
\]

Summing up the latter inequality and \rf{ai22}, we obtain
\bse
\|z\|_1&=& \sum_{i\in I \cap P_n} |z_i| +\sum_{i\in I_+ \cap P_+} z_i + \left[ \sum_{i\in I_- \cap P_+} { |z_i|} + \sum_{i \in \bar{I}} |z_i| \right] \le \nu + 2 \mu + 2 \kappa\\
&\le&
\nu + 2 \mu + \frac{2\beta(\e+\delta)+2\xi(\nu+(\theta+1)\mu)}{1-\xi} =  \frac{1+\xi}{1-\xi} \nu + \frac{2(1+\xi\theta)}{1-\xi}\mu + \frac{2\beta}{1-\xi} (\e+\delta),
\ese
which is \rf{eq256}.
\par
To show \rf{eq256'} observe that increasing $\e$ to $\e'=\e+\alpha\mu$, we can think that the true signal underlying the observation $y$ is $w^s$ rather than $w$; note that (\ref{Opt!}) implies that
\begin{equation}\label{Opt1!}
\|x\|_1\leq \|w^s\|_1+\nu',\,\,\nu'=\nu+\mu.
\end{equation}
We can now repeat the reasoning which follows (\ref{Opt!}), with (\ref{Opt1!}) in the role of (\ref{Opt!}), $w^s$ in the role of $w$, $\e'$ in the role of $\e$ and $0$ in the role of $\mu$, thus arriving at the following analogy of the bound (\ref{eq256}):
$$
\|{x}-w^s\|_1\leq \frac{1+\xi}{1-\xi} \nu' + \frac{2\beta}{1-\xi} (\e'+\delta),
$$
whence
$$
\|{x}-w\|_1\leq \frac{1+\xi}{1-\xi} \nu' + \frac{2\beta}{1-\xi} (\e'+\delta)+\mu,
$$
which is nothing but \rf{eq256'}. \qed

\section{Proof of Proposition \protect{\ref{propmain}}}\label{App:Propmain}

{\bf }
Let $A$ satisfy $\bVSG_s(\xi,\theta,\rho,\sigma)$, and let $Y=[y_1,...,y_n]$ and $v$ satisfy (\ref{condition}). Let, further, $I\subset\{1,...,n\}$ be such that $\Card(I)\leq s$, and let $x\in\bR^n$.
Let $u\in\bR^n$ be given by
$$
u_i=\left\{\begin{array}{ll}1+\theta\xi,&i\in P_+\cap I,~x_i\geq0\\
1-\xi,&i\in P_+\cap I,~x_i<0\\
(1+\xi)\sign(x_i),&i\in P_n\cap I\\
0,&i\not\in I\\
\end{array}\right..
$$
Note that $u$ has at most $s$ nonzero entries, the entries of $u$ with indices from $P_+$ belong to $[0,1+\theta\xi]$, and the modulae of entries in $u$ with indices from $P_n$ are $\leq1+\xi$, so that $u^Tz\leq\Phi_s(z)$ for all $z$. We have
\bse
\lefteqn{u^T[I-Y^TA]x = \sum\limits_iu^TC_i[Y,A]x_i= \sum\limits_{i:x_i\geq0}u^TC_i[Y,A]x_i+\sum\limits_{i:x_i<0}u^T[-C_i[Y,A]]|x_i|}\\
&\leq&\sum\limits_{i:x_i\geq0}\Phi_s(C_i[Y,A])x_i+\sum\limits_{i:x_i<0}\Phi_s(-C_i[Y,A])|x_i|\qquad\hbox{[since $u^Tz\leq \Phi_s(z)$]}\\
&\leq& \sum\limits_{i:x_i\geq0,i\not\in P_+}[\xi+(A^Tv)_i]x_i+\sum\limits_{i:x_i\geq0,i\in P_+}[\theta\xi+(A^Tv)_i]x_i+\sum\limits_{i:x_i<0}[\xi-(A^Tv)_i]|x_i|\quad
\hbox{[by (\ref{condition})]}\\
&=& \xi\left[\sum\limits_{i:x_i\geq0,i\not\in P_+}x_i+\theta\sum\limits_{i:x_i\geq0,i\in P_+}x_i+\sum\limits_{i:x_i<0}|x_i|\right]+x^TA^Tv\\
&=&\xi\left[\sum\limits_{i\in P_+}\max[-x_i,\theta x_i]+\sum\limits_{i\in P_n}|x_i|\right]+x^TA^Tv,
\ese
whence
\be
u^T[I-Y^TA]x \leq \xi\left[\sum\limits_{i\in P_+}\max[-x_i,\theta x_i]+\sum\limits_{i\in P_n}|x_i|\right]+\rho\|Ax\|
\ee{whence}
(recall that $\|v\|_*\leq\rho$).
On the other hand, recalling the definition of $u$ and that $\|y_i\|_*\leq\sigma$, we have
$$
\begin{array}{rcl}
u^T[I-Y^TA]x
&=&u^Tx-\sum\limits_{i\in I}u_iy_i^TAx
\\
&=& \sum\limits_{i\in I\cap P_+}\max[(1-\xi)x_i,(1+\theta\xi)x_i]+(1+\xi)\sum\limits_{i\in I\cap P_n}|x_i|
-\sum\limits_{i\in I}u_iy_i^TAx\\
&\geq&\sum\limits_{i\in I\cap P_+}\max[(1-\xi)x_i,(1+\theta\xi)x_i]+(1+\xi)\sum\limits_{i\in I\cap P_n}|x_i|\\
&&-\underbrace{\sigma\left[\sum\limits_{i\in I\cap P_+}(1+\theta\xi)
+\sum\limits_{i\in I\cap P_n}(1+\xi)\right]}_{\leq \beta-\rho}\|Ax\|.
\end{array}
$$
Combining the resulting inequality with (\ref{whence}), we get
$$
\sum\limits_{i\in I\cap P_+}\left[x_i+\xi\max[-x_i,\theta x_i]\right]+(1+\xi)\sum\limits_{i\in I\cap P_n}|x_i|\leq \beta\|Ax\|+
\xi\left[\sum\limits_{i\in P_+}\max[-x_i,\theta x_i]+\sum\limits_{i\in P_n}|x_i|\right]
$$
with $\beta$ given by (\ref{beta}), or, equivalently,
$$
\sum\limits_{i\in I\cap P_+}x_i +\sum\limits_{i\in I\cap P_n}|x_i|\leq \beta\|Ax\|+\xi\left[\sum\limits_{i\in P_+\setminus I}\max[-x_i,\theta x_i]+
\sum\limits_{i\in P_n\setminus I}|x_i|\right].
$$
The latter relation holds true for every $x\in\bR^n$ and for every set $I\subset\{1,...,n\}$ of cardinality $\leq s$, so that $A$
satisfies $\bSG_{s,\beta}(\xi,\theta)$. \qed

\section{Proof of Proposition \protect{\ref{limits}}}\label{App:Proplimits}
{\bf }
Proof is based on the following
\begin{lemma}\label{lemnew} Let $Z$ be a $\nu\times \nu$ matrix of rank $m$, $s>1$ be a positive integer, and $\delta_i\in(0,1]$, $1\leq i\leq \nu$, be such that
for the columns $C_i$ of the matrix $I_\nu-Z$ it holds $\|C_i\|_{s,1}\leq 1-\delta_i$. Assume that
\begin{equation}\label{besuch}
\nu>(2\sqrt{2m}+1)^2.
\end{equation}
Then
\begin{equation}\label{weneed100}
s\leq 2\sqrt{2m}+1.
\end{equation}
\end{lemma}
\def\Rank{{\hbox{\rm rank}}}
{\bf Proof of the lemma.} Let $\sigma_i=Z_{ii}$, and let $\gamma_i$ be the sum of $s-1$ largest magnitudes of the entries in $C_i$ with indices different from $i$. We have
$$
1-\sigma_i+\gamma_i\leq\|C_i\|_{s,1}\leq1-\delta_i,
$$
consequently $\sigma_i\geq\delta_i+\gamma_i>0$. Let us set $\lambda_i={1\over \sigma_i}$, and let $\bar{Z}$ be the matrix with the columns $\bar{Z}_i=\lambda_iZ_i$, where $Z_i$ is the $i$-th column in $Z$. Note that $\bar{Z}$ is of the same rank $m$ as $Z$, and that $\bar{Z}_{ii}=1$
for all $i$. Recalling that $\gamma_i<\sigma_i$, we have also
$$\|\bar{Z}_i\|_{s-1,1}=\lambda_i\|Z_i\|_{s-1,1}\leq\lambda_i[\gamma_i+\sigma_i]\leq 2\lambda_i\sigma_i=2. $$
Now let $\bar{s}=\min[s-1,\lfloor \nu^{1/2}\rfloor]$, so that $\bar{s}\geq1$ due to $s>1$. We have $\|\bar{Z}_i\|_{\bar{s},1}\leq\|\bar{Z}_i\|_{s-1,1}\leq2$ and $\bar{s}^2\leq \nu$. From the latter inequality and due to $\|\bar{Z}_i\|_2^2\leq \max\{1,\nu\bar{s}^{-2} \} \|\bar{Z}_i\|_{\bar{s},1}^2$ (cf. the proof of   \cite[Proposition 4.2]{JNCS}), it follows that  $\|\bar{Z}_i\|_2^2\leq 4\nu\bar{s}^{-2}$. We conclude that $\|\bar{Z}\|_2^2\leq 4\nu^2\bar{s}^{-2}$, where for a matrix $B$, $\|B\|_2$ is the Frobenius norm of $B$. Setting $H={1\over 2}[\bar{Z}+\bar{Z}^T]$, we have therefore $\|H\|_2^2\leq 4{\nu}^2\bar{s}^{-2}$. On the other hand, $\Tr(H)=\sum_{i=1}^\nu\bar{Z}_{ii}=\nu$, while $\Rank(H)\leq 2m$, whence, denoting by $\mu_i$, $1\leq i\leq p\leq 2m$, the nonzero eigenvalues of $H$, we have $$\|H\|_2^2=\sum_{i=1}^p\mu_i^2\geq
(\sum_{i=1}^p\mu_i)^2/p=(\Tr(H))^2/p\geq \nu^2/(2m).$$
We arrive at the inequality
$
4{\nu}^2\bar{s}^{-2}\geq\|H\|_2^2\geq \nu^2/(2m),
$
thereby

\begin{equation}\label{arriveat}
\bar{s}^2\leq8m.
\end{equation}
Assuming that $\bar{s}=\lfloor \nu^{1/2}\rfloor$, (\ref{arriveat}) says that ${\nu}\leq (2\sqrt{2m}+1)^2$, which is impossible. The only other option is that $\bar{s}=s-1$, and we arrive at (\ref{weneed100}). \qed
\paragraph{Lemma \ref{lemnew} $\Rightarrow$ Proposition \ref{limits}:} Let $Y,v$ satisfy (\ref{condition}). Consider first the case when $\nu:=\Card(P_n)\geq n/2.$ Denoting by $\widehat{C}_i$
the $\nu$-dimensional vector comprised of the last $\nu$ entries in $C_i=C_i[Y,A]$ (i.e., entries with indices from $P_n$). By (\ref{condition}), for every $i\in P_n$ and for every set $I\subset P_n$ with $\Card(I)\leq s$ we have
$$
\begin{array}{l}
\sum_{j\in I}(1+\xi)|[C_i]_j|\leq\Phi_s(-C_i)\leq \xi-(A^Tv)_i,\ \
\sum_{j\in I}(1+\xi)|[C_i]_j|\leq\Phi_s(C_i)\leq \xi+(A^Tv)_i,\\
\end{array}
$$
thus for any $i\in P_n$,
$$
2(1+\xi)\|\widehat{C}_i\|_{s,1}\leq \Phi_s(-C_i)+\Phi_s(C_i)\leq2\xi,
$$
so that $\|\widehat{C}_i\|_{s,1}< 1/2$. We see that the South-Eastern $\nu\times\nu$ submatrix $Z$ of $Y^TA$ satisfies the premise of Lemma \ref{lemnew}, while the size $\nu$ of $Z$ satisfies (\ref{besuch}) due to (\ref{eq119}) and $\nu\geq n/2$. Applying the lemma, we arrive at (\ref{eq191}).
\par
Now consider the case when $\Card(P_n)<n/2$, that is, $\nu:=\Card(P_+)\geq n/2$. By (\ref{condition}), setting $C_i=C_i[Y,A]$, for every set $I\subset P_+$ with $\Card(I)\leq s$ and every $i\in P_+$ we have
$$
\begin{array}{l}
\sum_{j\in I}(1+\theta\xi)\max[-[C_i]_j,0]\leq \Phi_s(-C_i)\leq \xi-(A^Tv)_i,\\
\sum_{j\in I}(1+\theta\xi)\max[[C_i]_j,0] \leq\Phi_s(C_i)\leq\theta\xi+(A^Tv)_i,\\
\end{array}
$$
whence
$$
\sum_{j\in I}|[C_i]_j|\leq {\xi(1+\theta)\over1+\theta\xi}<1.
$$
Since the latter inequality holds true for every subset $I$ of $P_+$ with $\Card(I)\leq s$, when denoting by $\bar{C}_i$ the part of $C_i$ comprised of the first $\nu$ entries (those with indexes from $P_+$), we have for all $i\in P_+$:
$$
\|\bar{C}_i\|_{s,1}<1.
$$
Now the proof can be completed exactly as in the previous case, with the North-Western $\nu\times \nu$ submatrix of $Y^TA$ in the role of $Z$. \qed
\section{Proof of Proposition \protect{\ref{Prop_NEMP}}}\label{App:Prop_NEMP}
{\bf }
Let us proceed by induction. First, let us show that $(a_{k-1},b_{k-1})$ implies $(a_{k},b_{k})$. Thus, assume that
$(a_{k-1},b_{k-1})$ holds true.
Let $z^{(k-1)}=w-v^{(k-1)}$. By $(a_{k-1})$, $z^{(k-1)}$ is supported on the support of $w$ and is such that $z^{(k-1)}_i\geq0$ for $i\in P_+$. Note that
\bse
z^{(k-1)}-u&=&w-v^{(k-1)}-Y^T(y-Av^{(k-1)})=(I-Y^TA)(w-v^{(k-1)})-Y^Te\\
&=&(I-Y^TA)z^{(k-1)}-Y^Te,
\ese
where $e=y-Aw$ with $\|Y^Te\|_\infty\le \sigma\delta$ due to (\ref{feq99}.c). Then by (\ref{feq99}.a,b) for any $i\in P_+$,
\bse
-\tau_-\left[\sum_{j\in P_+}z^{(k-1)}_j+\sum_{j\in P_n}|z^{(k-1)}_j|\right]-\sigma\delta\le z^{(k-1)}_i -u_i\le
\tau_+\left[\sum_{j\in P_+}z^{(k-1)}_j + \sum_{j\in P_n}|z^{(k-1)}_j|\right]+\sigma\delta,
\ese
consequently,
\be
 -\gamma_-:=-\tau_- \alpha_{k-1}-\sigma\delta\leq z^{(k-1)}_i-u_i\leq \gamma_+:=\tau_+\alpha_{k-1}+\sigma\delta.
\ee{cover}
We conclude that for any $i\in P_+$ the interval $S_i=[u_i-\gamma_-,\,u_i+\gamma_+]$ of the width
$$
\ell_+=[\tau_-+\tau_+]\alpha_{k-1}+2\sigma\delta,
$$
covers $z^{(k-1)}_i$. In the same way for any $i\in P_n$
\[
-\gamma:=-\tau\alpha_{k-1}-\sigma\delta\le z^{(k-1)}_i -u_i\le
\tau\alpha_{k-1}+\sigma\delta=\gamma,
\]
so that the interval $S_i=[u_i-\gamma,\,u_i+\gamma]$ of the width
$$
\ell=2\tau\alpha_{k-1}+2\sigma\delta,
$$
covers $z^{(k-1)}_i$ when $i\in P_n$.

Recalling that $z_i^{(k-1)}\geq0$ for $i\in P_+$, %i.e., $u_i+\gamma_i\geq0$ for $i\in P_+$,
the closest to $0$ point of $S_i$ is
\bse
\begin{array}{rlllcrlll}
\widetilde{\Delta}_i&=&[u_i-\gamma_-]_+&\mbox{for }i\in P_+,&&
\widetilde{\Delta}_i&=&[u_i-\gamma]_+&\mbox{for }i\in P_n,\;\;u_i\ge 0,\\
\widetilde{\Delta}_i&=&-[|u_i|-\gamma]_+&\mbox{for }i\in P_n,\;\;u_i<0,
\end{array}
\ese
that is, $\widetilde{\Delta}_i=\Delta_i$ for all $i$. Since the segment $S_i$ covers $z^{(k-1)}_i$ and $\Delta_i$ is the closest to 0 point in $S_i$, while the width of $S_i$ is at most $\ell\vee\ell_+$, we clearly have
\be\begin{array}{rlcrl}
(a)&\Delta_i\in \Conv\left\{0,z^{(k-1)}_i\right\},&&
(b)&|z^{(k-1)}_i-\Delta_i|\le \ell\vee\ell_+.\\
\end{array}
\ee{almost}
Since $(a_{k-1})$ is valid, (\ref{almost}.a) implies that
\[
v^{(k)}_i=v^{(k-1)}_i+\Delta_i\in\left[ v^{(k-1)}_i+\Conv\left\{0,w-v^{(k-1)}_i\right\}\right]\subseteq\Conv\{0,w_i\},
\]
and $(a_k)$ holds.
Further, let $I$ be the support of $w^s$. Relation $(a_k)$ clearly implies that $|z^{(k)}_i|\le |w_i|$, and we can write due to (\ref{almost}.b):
\bse
\|w-v^{(k)}\|_1&=&\sum_{i\in I} |w-[v^{(k-1)}_i+\Delta_i]|+\sum_{i\not\in I}|z^{(k)}_i| \\
&\le& \sum_{i\in I} |z^{(k-1)}_i-\Delta_i|+\sum_{i\not\in I}|w_i|\le s[\ell\vee\ell_+]+\mu=\alpha_k,
\ese
 which is $(b_k)$. The induction step is justified.
\par
It remains to show that $(a_0,b_0)$ holds true. Since $(a_0)$ is evident, all we need is to justify $(b_0)$.  Let
\[
\alpha_*=\|w\|_1,
\]
and let $u=Y^Ty$. Same as above (cf. \rf{cover}), we
have for all $i$:
$$
|w_i-u_i|\le \max\{\tau_-,\tau_+,\tau\}\alpha_*+\sigma\delta={\rho\over s}\alpha_*+\sigma\delta.
$$
Then
\bse
\alpha_*=\sum_{i\in I}|w_i|+\sum_{i\not \in I}|w_i|\leq \sum_{i\in I} [|u_i|+
{\rho\over s}\alpha_*+\sigma\delta]+\mu\le \|u\|_{s,1}+
\rho\alpha_*+s\sigma\delta+\mu.
\ese
Hence
\[
\alpha_*\leq \alpha_0={\|u\|_{s,1}+s\sigma\delta+\mu\over 1-\rho},
\]
which implies $(b_0)$.   \qed

\newpage
\section{ONLINE SUPPLEMENT}

\subsection{Proof of Proposition \protect{\ref{propnew}}}\label{App:Monotonicity}

Let $Y=[Y_1,...,Y_n],v,\sigma,\rho$ certify the validity of $\bVSG^*_{s,\beta}(\xi,\theta)$, and let $\beta'\geq\beta$, $\theta'\geq\theta$ and $\xi'\in[\xi,1)$. Let us set
$$
\lambda={1+\theta\xi\over 1+\theta'\xi'},\,\,\mu={1+\xi\over 1+\xi'}.
$$
so that $\lambda,\mu\in[0,1]$, and let $Y'$ be as in the assertion to be proved, that is, the columns of $Y'$ are multiples of those of $Y$: $Y^\prime_i=\lambda Y_i$ when $i\in P_+$ and $Y^\prime_i=\mu Y_i$ otherwise. All we need to prove is that $(Y',v,\sigma,\rho)$ certify the validity of $\bVSG^*_{s,\beta'}(\xi',\theta')$, and this immediately reduces to verification of the following fact:
\begin{lemma}\label{lemlast}
Let $i$, $1\leq i\leq n$, be fixed, and let $z\in\bR^n $  for any $I\subset \{1,...,n\}$ of cardinality $s$ satisfy the relations
\begin{equation}\label{given}
\begin{array}{ll}
(a)&(1+\theta\xi)\sum\limits_{j\in P_+\cap I}\max[z_j-\delta_{ij},0] +(1+\xi)\sum\limits_{j\in P_n\cap I}|z_j-\delta_{ij}|+(Av)_i\leq \xi,\\
(b)&(1+\theta\xi)\sum\limits_{j\in P_+\cap I}\max[\delta_{ij}-z_j,0] +(1+\xi)\sum\limits_{j\in P_n\cap I}|z_j-\delta_{ij}|-(Av)_i\\
&\leq \eta=\left\{\begin{array}{ll}\theta\xi,&i\in P_+,\\
\xi,&i\in P_n,\\
\end{array}\right.\\
\end{array}
\end{equation}
where $\delta_{ij}=\left\{\begin{array}{ll}0,&j\neq i,\\1,&i=j.\end{array}\right.$ %Let us set $w_j=\lambda z_j$, $j\in P_+$, and $w_j= \mu z_j$, $j\in P_n$.
Then
%$$
%\begin{array}{ll}
%\multicolumn{2}{l}{\forall I\subset \{1,...,n\},\Card(I)=s:}\\
%(a)&(1+\theta'\xi')\sum\limits_{j\in P_+\cap I}\max[w_j-\delta_{ij},0] +(1+\xi')\sum\limits_{j\in P_n\cap I}|w_j-\delta_{ij}|+(Av)_i\leq \xi'\\
%(b)&(1+\theta'\xi')\sum\limits_{j\in P_+\cap I}\max[\delta_{ij}-w_j,0] +(1+\xi')\sum\limits_{j\in P_n\cap I}|w_j-\delta_{ij}|-(Av)_i\\
%&\leq \eta_+=\left\{\begin{array}{ll}\theta'\xi',&i\in P_+\\
%\xi',&i\in P_n\\
%\end{array}\right.\\
%\end{array}
%$$
%or, which is the same,
for every set $I\subset\{1,...,n\}$ of cardinality $s$ we have
\begin{equation}\label{toprove}
\begin{array}{ll}
(a)&(1+\theta'\xi')\sum\limits_{j\in P_+\cap I}\max[\lambda z_j-\delta_{ij},0] +(1+\xi')\sum\limits_{j\in P_n\cap I}|\mu z_j-\delta_{ij}|+(Av)_i\leq \xi',\\
(b)&(1+\theta'\xi')\sum\limits_{j\in P_+\cap I}\max[\delta_{ij}-\lambda z_j,0] +(1+\xi')\sum\limits_{j\in P_n\cap I}|\mu z_j-\delta_{ij}|-(Av)_i\\
&\leq \eta_+=\left\{\begin{array}{ll}\theta'\xi',&i\in P_+,\\
\xi',&i\in P_n.\\
\end{array}\right.\\
\end{array}
\end{equation}
\end{lemma}
{\bf Proof.}
Taking into account the definition of $\lambda,\mu$, in the case of $i\not\in I$ the relations (\ref{toprove}) are readily given by (\ref{given}), hence we can assume $i\in I$.
Consider two possible cases: $i\in P_+\cap I$ and $i\in P_n\cap I$.
\paragraph{The case of $i\in P_+\cap I$.} In  this case (\ref{given}) reads:
\begin{equation}\label{given1}
\begin{array}{ll}
(a)&(1+\theta\xi)\max[z_i-1,0]+(1+\theta\xi)\sum\limits_{j\in P_+\cap I,j\neq i}\max[z_j,0] \\ &+(1+\xi)\sum\limits_{j\in P_n\cap I}|z_j|+(Av)_i\leq \xi,\\
(b)&(1+\theta\xi)\max[1-z_i,0]+(1+\theta\xi)\sum\limits_{j\in P_+\cap I,j\neq i}\max[-z_j,0] \\ &+(1+\xi)\sum\limits_{j\in P_n\cap I}|z_j|-(Av)_i\leq \theta\xi,
\end{array}
\end{equation}
and our goal is to verify that then
\begin{equation}\label{toprove1}
\begin{array}{ll}
(a)&(1+\theta'\xi')\max[\lambda z_i-1,0]\\&+\overbrace{(1+\theta'\xi')\lambda}^{=1+\theta\xi}\sum\limits_{j\in P_+\cap I,j\neq i}\max[z_j,0]
+\overbrace{(1+\xi')\mu}^{=1+\xi}\sum\limits_{j\in P_n\cap I}|z_j|+(Av)_i\leq \xi',\\
(b)&(1+\theta'\xi')\max[1-\lambda z_i,0]\\&+\underbrace{(1+\theta\xi)\sum\limits_{j\in P_+\cap I,j\neq i}\max[-z_j,0]+(1+\xi)\sum\limits_{j\in P_n\cap I}|z_j|-(Av)_i}_{:=R}\leq \theta'\xi'.
\end{array}
\end{equation}
We have $\lambda z_i-1\leq\lambda(z_i-1)$ due to $\lambda\leq1$, consequently \[
\max[\lambda z_i-1,0]\leq \max[\lambda(z_i-1),0]=\lambda\max[z_i-1,0], \]
and therefore (\ref{toprove1}.$a$) follows from (\ref{given1}.$a$) due to
$(1+\theta'\xi')\lambda=1+\theta\xi$ and $\xi'\geq\xi$. It remains to verify (\ref{toprove1}.$b$). Assume, first, that $\lambda z_i\leq1$. From (\ref{given1}.$b$) it follows that
\[
(1+\theta\xi)[1-z_i]+R\leq (1+\theta\xi)\max[1-z_i,0]+R\leq \theta\xi,\] implying
$
z_i\geq {1+R\over 1+\theta\xi}$ and therefore \[
1-\lambda z_i\leq 1-{1+R\over1+\theta'\xi'}={\theta'\xi'-R\over 1+\theta'\xi'}.\]
 Since we are in the case $1-\lambda z_i\geq0$, we arrive at
$$
(1+\theta'\xi')\max[1-\lambda z_i,0]+R=(1+\theta'\xi')[1-\lambda z_i]+R\leq (1+\theta'\xi'){\theta'\xi'-R\over 1+\theta'\xi'}+R=\theta'\xi',
$$
as required in (\ref{toprove1}.$b$). The case of $1-\lambda z_i\leq0$ is trivial, since here the left hand side in (\ref{toprove1}.$b$) clearly is $\leq$ the left hand side in (\ref{given1}.$b$),
while $\theta'\xi'\geq\theta\xi$, so that (\ref{toprove1}.$b$) is readily given by (\ref{given1}.$b$). Thus, when $i\in P_+\cap I$, (\ref{toprove1}) follows from (\ref{given1}).
\paragraph{The case of $i\in P_n\cap I$.} In this case (\ref{given}) means that
\begin{equation}\label{given2}
\begin{array}{ll}
(a)&(1+\theta\xi)\sum\limits_{j\in P_+\cap I,j\neq i}\max[z_j,0] +(1+\xi)|1-z_i|+(1+\xi)\sum\limits_{j\in P_n\cap I,j\neq i}|z_j|+(Av)_i\leq \xi,\\
(b)&(1+\theta\xi)\sum\limits_{j\in P_+\cap I}\max[-z_j,0] +(1+\xi)|1-z_i|+(1+\xi)\sum\limits_{j\in P_n\cap I,j\neq i}|z_j|-(Av)_i\leq \xi,\\
\end{array}
\end{equation}
and our goal is to verify that then
\begin{equation}\label{toprove2}
\begin{array}{ll}
(a)&(1+\theta'\xi')\sum\limits_{j\in P_+\cap I,j\neq i}\max[\lambda z_j,0] \\&
+(1+\xi')|1-\mu z_i|+(1+\xi')\mu\sum\limits_{j\in P_n\cap I,j\neq i}|z_j|+(Av)_i\leq \xi',\\
(b)&(1+\theta'\xi')\sum\limits_{j\in P_+\cap I}\max[-\lambda z_j,0] \\&
+(1+\xi')|1-\mu z_i|+(1+\xi')\sum\limits_{j\in P_n\cap I,j\neq i}|\mu z_j|-(Av)_i\leq \xi'.
\end{array}
\end{equation}
Comparing (\ref{given2}.$a$) with (\ref{toprove2}.$a$), and (\ref{given2}.$b$) with (\ref{toprove2}.$b$), we see that all we need in order to derive (\ref{toprove2}) from (\ref{given2})
is to verify the following statement: if $(1+\xi)|1-z|\leq\xi+a$, then $(1+\xi')|1-\mu z|\leq \xi'+a$. This is immediate: assuming $(1+\xi)|1-z|\leq\xi+a$, the premises in the following two implication chains hold true:
$$
\begin{array}{l}
(1+\xi)[1-z]\leq\xi+a\Rightarrow z\geq {1-a\over 1+\xi}\Rightarrow \mu z\geq {1-a\over 1+\xi'}\Rightarrow 1-\mu z\leq 1-{1-a\over 1+\xi'}={\xi'+a\over 1+\xi'}\\
\qquad\qquad\Rightarrow (1+\xi')[1-\mu z]\leq \xi'+a,\\
(1+\xi)[z-1]\leq\xi+a\Rightarrow z\leq 1+{\xi+a\over 1+\xi}\Rightarrow \mu z\leq {1+2\xi+a\over 1+\xi'}\Rightarrow \mu z-1\leq {2\xi-\xi'+a\over 1+\xi'}\\
\qquad\qquad\Rightarrow
(1+\xi')[\mu z-1] \leq 2\xi-\xi'+a\Rightarrow (1+\xi')[\mu z-1]\leq \xi'+a,\\
\end{array}
$$
while the resulting inequalities in these chains lead to the desired conclusion $(1+\xi')|1-\mu z|\leq \xi'+a$. \qed %
\subsection{``Trigonometric polynomials'' example}\label{LemmaMoved}
The validity of the claim concluding Section \ref{sec:limits} is readily given by the following
\begin{lemma}\label{lemma:cyclic}
For any positive integer $d$, let $n\geq 4\pi d$, and $A$ be the matrix obtained from the basic trigonometric polynomials as described in Section \ref{sec:limits}, then the condition $\bVSG_s(\xi,\theta,\rho,\sigma)$ can hold true for $s\leq2$ only.
\end{lemma}
{\bf Proof.} Let $L$ be the $n\times n$ permutation matrix corresponding to the cyclic shift $e_j\mapsto e_{j_+}$, $j_+=(j+1)\,\hbox{mod}\,n$, of the standard basic orths $e_0,...,e_{n-1}$ in $\bR^n$, and $R$ be the $m\times m$ orthogonal block-diagonal matrix with the North-Western block $1$ and $d$ additional $2\times 2$ diagonal blocks $\left[\begin{array}{cc}\cos(2\pi i/n)&
-\sin(2\pi i/n)\cr\sin(2\pi i/n)&\cos(2\pi i/n)\end{array}\right]$, $1\leq i\leq d$. Denoting by $A_j$ the $j$-th column of $A$, $0\leq j\leq n-1$, we clearly have $R A_j=A_{j_+}$, hence $A=RAL^{-1}$ and therefore also $A=R^iAL^{-i}$ for $1\leq i\leq n$. Now assume that $Y,v$ satisfy (\ref{condition}) for certain $\xi<1$, $\theta\ge1$, $\rho$, $\sigma$. Then
$$
\max\limits_i\left[\Phi_s(-C_i[Y,A])+\Phi_s(C_i[Y,A])\right]\leq \xi(1+\theta),
$$
in this way, it is immediately seen, $\max_i\|C_i[Y,A]\|_{s,1}\leq \kappa:={\xi(1+\theta)\over 1+\theta\xi}<1$, or, which is the same,
\[
\Gamma(I-Y^TA)\leq\kappa<1,
\]
where $\Gamma(Z)$ is the maximum of the $\|\cdot\|_{s,1}$-norms of columns of $Z\in\bR^{n\times n}$. Observe that $\Gamma$ is a convex function which is symmetric in the sense that
$
\Gamma(PZP^T)=\Gamma(Z)
$
whenever $P$ is a permutation matrix. Now let $\bar{Y}={1\over n}\sum_{i=1}^n R^{-i}YL^i$. Since $L^n=I_n$, $R^{-n}=I_m$, we have $R^{-1}\bar{Y}L=\bar{Y}$. We claim that
\[
\Gamma(I-\bar{Y}^TA)\leq\kappa.
\]
Indeed, we have
\bse
\Gamma(I-\bar{Y}^TA)&=&\Gamma({1\over n}\sum_{i=1}^n [I-L^{-i}Y^TR^iA])\\
&\leq& {1\over n}\sum_{i=1}^n\Gamma(I-L^{-i}Y^TR^iA)\;\;\;\;\hbox{\ [since $\Gamma$ is convex]}\\
&=&{1\over n}\sum_{i=1}^n\Gamma(L^{-i}\left[I-Y^T[R^iAL^{-i}]\right]L^i)\\
&=&{1\over n}\sum_{i=1}^n\Gamma(I-Y^TA)\;\;\;\;\hbox{\ [since $\Gamma$ is symmetric and $R^iAL^{-i}=A$]} \\
&=&\Gamma(I-Y^TA)\ese
Now let
\[
y_j(\phi)=\bar{Y}_{0j}+\sum_{i=1}^d[\bar{Y}_{2i-1,j}\cos(i\phi)+\bar{Y}_{2i,j}\sin(i\phi)].
\]
We have $R^{-1}\bar{Y}L=\bar{Y}$, that is, $R^{-1}\bar{Y}=\bar{Y}L^{-1}$. In other words, the columns $\bar{Y}_j$ of $\bar{Y}$ satisfy the relation $\bar{Y}_j=R\bar{Y}_{j_-}$, where $j_-=(j-1)\,\hbox{mod}\,n$. This is nothing but $y_j(\phi)\equiv y_{j_-}(\phi-\delta)$, $\delta=2\pi/n$, whence $y_j(\phi)=y_0(\phi-j\delta)$. Observe that the $j$-th column in $\bar{Y}^TA$ has the entries
\[
\bar{Y}_i^TA_j=y_i(j\delta)=y_0((j-i)\delta), \;\;0\leq i\leq n-1,
 \]
meaning that the columns in the matrix $I-\bar{Y}^TA$ are cyclic shifts of each other (so that the
\hbox{$\|\cdot\|_{s,1}$}-norms of all columns are the same), and the zero column is comprised of the values of the trigonometric polynomial $1-y_0(\phi)$ on the grid $G=\{\phi_j={2\pi j\over n}:0\leq j<n\}$. Assuming $s>1$, when denoting by $\gamma$ the sum of $s-1$ largest magnitudes of entries in the $(n-1)$-dimensional vector $\{y_0(\phi_i)\}_{i=1}^{n-1}$, we have
\[
1-y_0(0)+\gamma\leq \|C_0[\bar{Y},A]\|_{s,1}\leq\kappa<1,
\]
thereby $\mu:=y_0(0)>\gamma$. Now let $M=\max\limits_{0\leq\phi\leq 2\pi}|y_0(\phi)|$, and let $\bar{\phi}\in\Argmax_\phi|y_0(\phi)|$, so that $y_0^\prime(\bar{\phi})=0$. By Bernstein theorem, we have $|y^{\prime\prime}_0(\phi)|\leq d^2M$ for all $\phi$, whence $|y_0(\phi)|\geq M/2$ when $|\phi-\bar{\phi}|\leq 1/d$, so that
\[
\Card\{j:|y_0(\phi_j)|\geq M/2\}>{n\over \pi d}-1.
\]
It follows that $\gamma\geq \min\left[s-1,{n\over \pi d}-2\right]M/2$, while $\mu=y_0(0)\leq M$. Thus, the relation $\mu>\gamma$ implies that
$$
\min[s-1,{n\over \pi d}-2]<2,
$$
that is, $s\leq 2$ provided that $n\geq 4\pi d$. \qed

\end{document}